\documentclass[11pt,reqno,twoside]{article}
\usepackage{amsmath,amssymb,theorem}
\theorembodyfont{\itshape}
\newtheorem{thm}{Theorem}[section]
\newtheorem{lem}[thm]{Lemma}
\newtheorem{prop}[thm]{Proposition}
{\theorembodyfont{\upshape}

\newtheorem{rem}[thm]{Remark}
}

\newcommand{\Proof}[1][]{\noindent{\itshape Proof#1. }}
\newcommand{\EndProof}{~$\Box$\bigskip}
\makeatletter
    
    \@addtoreset{equation}{section}
\makeatother

\def\mbN{{\mathbb N}}
\def\mbR{{\mathbb R}}
\def\mbC{{\mathbb C}}

\def\mbZ{{\mathbb Z}}

\def\mcE{\mathcal{E}}
\def\mcF{\mathcal{F}}

\def\mcH{\mathcal{H}}

\def\mcX{\mathcal{X}}

\def\ua{\uparrow}

\def\ol{\overline}

\def\sm{\setminus}
\def\sbs{\subset}

\def\sps{\supset}
\def\ssbs{\Subset}
\def\ptl{\partial}
\def\es{\emptyset}

\def\wt{\widetilde}

\def\a{\alpha}    \def\b{\beta}              
\def\D{\Delta}    \def\e{\varepsilon}        \def\Ph{\Phi}
  \def\g{\gamma}           
   \def\L{\Lambda}     \def\m{\mu}        \def\n{\nu}
               \def\O{\Omega}
              \def\s{\sigma}     
          
\def\x{\xi}       
\def\z{\zeta}           

\begin{document}
\title{A Variational Principle in the Dual Pair of Reproducing Kernel Hilbert Spaces and an Application}
\author{ Hyun Jae Yoo\footnote{University College, Yonsei University, 134
Shinchon-dong, Seodaemoon-gu, Seoul 120-749, Korea. E-mail:
yoohj@yonsei.ac.kr} }
\date{}
  \maketitle

\begin{abstract}
Given a positive definite, bounded linear operator $A$ on the
Hilbert space $\mcH_0:=l^2(E)$, we consider a reproducing kernel
Hilbert space $\mcH_+$ with a reproducing kernel $A(x,y)$. Here
$E$ is any countable set and $A(x,y)$, $x,y\in E$, is the
representation of $A$ w.r.t. the usual basis of $\mcH_0$. Imposing
further conditions on the operator $A$, we also consider another
reproducing kernel Hilbert space $\mcH_-$ with a kernel function
$B(x,y)$, which is the representation of the inverse of $A$ in a
sense, so that $\mcH_-\sps\mcH_0\sps\mcH_+$ becomes a rigged
Hilbert space. We investigate a relationship between the ratios of
determinants of some partial matrices related to $A$ and $B$ and
the suitable projections in $\mcH_-$ and $\mcH_+$. We also get a
variational principle on the limit ratios of these values. We
apply this relation to show the Gibbsianness of the determinantal
point process (or fermion point process) defined by the operator
$A(I+A)^{-1}$ on the set $E$. It turns out that the class of
determinantal point processes that can be recognized as Gibbs
measures for suitable interactions is much bigger than that
obtained by Shirai and Takahashi.
\end{abstract}
{\bf Running head}. {A Variational Principle in RKHS's}


\section{Introduction}\label{sec:introduction}
In this paper we will consider certain variational principle
arising in the dual pair of reproducing kernel Hilbert spaces
(abbreviated RKHS's hereafter). Then we will find an application
in showing the Gibbsianness of some determinantal point processes
(in short DPP's) in discrete spaces.

Let $E$ be any countable set, e.g., $E=\mbZ^d$, the $d$-dimensional lattice space. Let $\mcH_0:=l^2(E)$ be the
space of square summable functions (sequences) on $E$ with inner product
\begin{equation}
(f,g)_0:=\sum_{x\in E}\overline{f(x)}g(x).
\end{equation}
We denote the corresponding norm by $\|\cdot\|_0$. Let $A$ be a
bounded positive definite operator on $\mcH_0$. We assume that the
kernel space is trivial: $\text{ker}A=\{0\}$. Then the range,
$\text{ran}A$, is dense in $\mcH_0$ and we introduce two new
norms. First on $\mcH_0$ we define
\begin{equation}\label{eq:-norm}
\|f\|_-^2:=(f,Af)_0,\quad f\in \mcH_0.
\end{equation}
Let $\mcH_-$ be the closure of $\mcH_0$ w.r.t. this norm. Next we define a norm $\|\cdot\|_+$ on $\text{ran}A$
by
\begin{equation}\label{eq:+norm}
\|g\|_+^2:=(g,A^{-1}g)_0,\quad g\in \text{ran}A.
\end{equation}
The closure of $\text{ran}A$ w.r.t. the norm $\|\cdot\|_+$ is
denoted by $\mcH_+$. We then get a triple with inclusions:
\begin{equation}\label{eq:rigging}
\mcH_-\sps\mcH_0\sps\mcH_+.
\end{equation}

Let us denote by ${\sf B}:=\{e_x:\,x\in E\}$ the usual basis of $\mcH_0$, i.e., $e_x\in \mcH_0$ is the unit
vector whose component is one at $x$ and zero at all other sites. Let $A(x,y)$, $x,y\in E$, be the
representation of $A$ w.r.t. the basis ${\sf B}$. Then the space $\mcH_+$ is nothing but a RKHS with a
reproducing kernel (abbreviated RK) $A(x,y)$, $x,y\in E$ (see Subsection 2.1). We allow $0\in \text{spec}A$, the
spectrum of $A$. That is, the inverse $A^{-1}$ of $A$ may be an unbounded operator on $\mcH_0$. But we will
impose some conditions on $A$ so that the space $\mcH_-$ is also a RKHS with a RK $B(x,y)$, $x,y\in E$. See the
hypothesis (H) in Section 2. Informally saying, the function $B(x,y)$ is the kernel function of the inverse
operator $A^{-1}$:
\begin{equation}\label{eq:kernel_B}
B(x,y)=A^{-1}(x,y),\quad x,y\in E.
\end{equation}

The variational principle we will address is the following. We notice first that the assumption of $\mcH_-$
being a RKHS implies in particular that $e_x\in \mcH_+$ for all $x\in E$ (see Subsection 2.1). Let $x_0\in E$ be
a fixed point and let $\{x_0\}\cup R_1\cup R_2=E$ be any partition of $E$. We define
\begin{eqnarray}\label{eq:alpha}
\a&:=&\lim_{\L\ua E}\a_\L;\nonumber\\
\a_\L&:=&\inf_{f\in \text{span}\{e_x:x\in \L\cap
R_1\}}\|e_{x_0}-f\|_-^2,
\end{eqnarray}
and similary
\begin{eqnarray}\label{eq:beta}
\b&:=&\lim_{\L\ua E}\b_\L;\nonumber\\
\b_\L&:=&\inf_{g\in \text{span}\{e_x:x\in \L\cap R_2\}}\|e_{x_0}-g\|_+^2,
\end{eqnarray}
where $\L$ increases to $E$ through finite subsets. We will show that the two numbers $\a$ and $\b$ are the
inverses to each other (Theorem \ref{thm:variational_principle}):
\begin{equation}\label{eq:variation_principle}
\a\b=1.
\end{equation}
This result has been shown by Shirai and Takahashi \cite{ST2} in
the case when $A$ is a strictly positive operator, and hence
$A^{-1}$ is also bounded. They applied this result to show the
Gibbsianness of a DPP defined by the operator $A(I+A)^{-1}$ (see
Section 2 for the definition of DPP's). In fact, the variational
principle \eqref{eq:variation_principle} will guarantee the
existence of global Papangelou intensity. In other words, it will
prove the existence of the limit of local Papangelou intensities
as the local region increases to the whole space (Theorem
\ref{thm:global_PI}). This proves the Gibbsianness of the DPP and
we will give a proper interaction potential and show also the
uniqueness of the Gibbs measure (Theorem \ref{thm:Gibbs_measure}).
The interaction potential is actually given by the logarithm of
the determinants of the submatrices of $A$:
\begin{equation}\label{eq:interaction_potential}
V(\x)=-\log \det (A(x_i,x_j))_{i,j=1}^n,
\end{equation}
where $\x=\{x_1,\cdots,x_n\}\sbs E$ is any finite configuration.

We remark here that the main idea in showing the Gibbsianness has
been borrowed from \cite{ST2}. We should, however, point out that
since the operators $A$ dealt with in \cite{ST2} are strictly
positive, there is a severe restriction in applications. For
example, if $A$ is a diagonal matrix with diagonal elements
$\a_x>0$ that decrease to zero as $x\to \infty$ (we let $E=\mbZ$
or $\mbN$), then the DPP corresponding to the operator
$A(I+A)^{-1}$ is clearly a Gibbs measure. The system has the
one-body interactions only and the potential energy is given by
\begin{equation}\label{eq:interaction_potential_diagonal}
V(\x)=-\sum_{x\in \x}\log \a_x.
\end{equation}
Even this kind of simple example lies outside the regime of \cite{ST2}. This paper improves \cite{ST2} (in
regard of Gibbsianness of DPP's) in that our setting includes more general classes as well as the above example.

This paper is organized as follows. In Section 2, we introduce the
basics of the RKHS's (Subsection 2.1) and DPP's (Subsection 2.2),
and then give the main results (Subsection 2.3). Section 3 is
devoted to the proof of variational principle, Theorem
\ref{thm:variational_principle}. In Section \ref{sec:proof of
global_PI}, we first prove the existence of the global Papangelou
intensity, Theorem \ref{thm:global_PI}. Then we prove the
Gibbsianness and its uniqueness, Theorem \ref{thm:Gibbs_measure}.
In the Appendix, we provide with some examples.

\section{Preliminaries and Main Results}
In this Section we review some basics of RKHS's and DPP's. Then we state the main results of this paper.

\subsection{Reproducing Kernel Hilbert Spaces}
For our convenience, we start from a Hermitian positive definite bounded linear operator $A$ on the complex
Hilbert space $\mcH_0:=l^2(E)$ equipped with an inner product
\begin{equation}\label{eq:0inner_product}
(f,g)_0:=\sum_{x\in E}\overline{f(x)}g(x),\quad f,g\in \mcH_0.
\end{equation}
Here $E$ is any countable set. Throughout this paper we assume
that the kernel space of $A$ is trivial:
\begin{equation}\label{eq:null_kernel}
\text{ker}A=\{0\}.
\end{equation}
Then, since $\overline{\text{ran}A}=(\text{ker}A^*)^\bot=(\text{ker}A)^\bot=\mcH_0$, $\text{ran}A$ is dense in
$\mcH_0$. As in the introduction, let ${\sf B}=\{e_x:\, x\in E\}$ be the usual basis of $\mcH_0$. Let $A(x,y)$,
$x,y\in E$, be the matrix element of the operator $A$ w.r.t. the basis ${\sf B}$:
\begin{equation}\label{eq:A_matrix}
A(x,y):=(e_x,Ae_y)_0,\quad x,y\in E.
\end{equation}
On the dense subspace $\text{ran}A$, we define a new inner product as
\begin{equation}\label{eq:+inner_product}
(f,g)_+:=(f,A^{-1}g)_0,\quad f,g\in \text{ran}A.
\end{equation}
Denote by $\|\cdot\|_+$ the resulting norm and let $\mcH_+$ be the
completion of $\text{ran}A$ w.r.t. $\|\cdot\|_+$. We notice that
$\mcH_+$ is a RKHS \cite{Ar, HY, HKPS} with kernel function
$A(x,y)$, that is the following defining conditions are satisfied:
\begin{itemize}
\item[(i)] For every $x\in E$, the function $A(\cdot,x)$ belongs
to $\mcH_+$,
 \item[(ii)] The reproducing property: for every $x\in E$ and $g\in \mcH_+$,
 \begin{equation}
 g(x)=(A(\cdot,x),g)_+.
 \end{equation}
\end{itemize}

Let us now consider another Hilbert space $\mcH_-$ which is the
closure of $\mcH_0$ w.r.t. the norm $\|\cdot\|_-$ induced by the
inner product:
\begin{equation}\label{eq:-inner_product}
(f,g)_-:=(f,Ag)_0,\quad f,g\in \mcH_0.
\end{equation}
It is important to notice that though $\mcH_0$ may be understood
as a class of functions defined on the set $E$, the completed
space $\mcH_-$ may not be a space of functions defined on the same
space $E$. This is so called a functional completion problem
\cite{Ar} and will be discussed below. By the boundedness of $A$
we have the inclusions:
\begin{equation}\label{eq:inclusions}
\mcH_-\sps\mcH_0\sps\mcH_+.
\end{equation}
We want to see $\mcH_-$ also as a RKHS. First we define a dual pairing between the spaces $\mcH_-$ and $\mcH_+$.
For $f\in \mcH_0$ and $g\in \text{ran}A$, define
\begin{equation}\label{eq:dual_pair}
_-\langle f,g\rangle_+:=\sum_{x\in E}\ol{f(x)}g(x).
\end{equation}
We have then the bound $|_-\langle  f,g\rangle_+|\le
\|f\|_-\|g\|_+$. Since $\mcH_0$ and $\text{ran}A$ are dense
respectively in $\mcH_-$ and $\mcH_+$, the dual pairing extends
continuously to a bilinear form on $\mcH_-\times \mcH_+$, for
which we use the same notation $_-\langle  f,g\rangle_+$, $f\in
\mcH_-$ and $g\in \mcH_+$, and the bound also continues to hold:
\begin{equation}\label{eq:bound_of_dual_pair}
|_-\langle  f,g\rangle_+|\le \|f\|_-\|g\|_+,\quad f\in
\mcH_-,\,\,g\in \mcH_+.
\end{equation}
For a convenience, we also define its conjugate bilinear form
\begin{equation}\label{eq:conjugate_form}
_+\langle  g,f\rangle_-:=\ol{_-\langle  f,g\rangle_+},\quad f\in
\mcH_-,\,\,g\in \mcH_+.
\end{equation}
Notice that for $f\in \mcH_0$, $Af\in \mcH_+$ and
\begin{equation}\label{eq:bound_of_A}
\|Af\|_+^2=(Af,A^{-1}Af)_0=\|f\|_-^2.
\end{equation}
Thus, $A$ extends to an isometry between $\mcH_-$ and $\mcH_+$. We
will denote the extension by the same $A$ and its inverse by
$A^{-1}$.

Let us now introduce the notion of {\it functional completion} of
an incomplete class   $\sf F$ of functions on $E$ which is a
pre-Hilbert space. By this, as introduced in \cite[p 347]{Ar}, we
mean a completion of $\sf F$ by adjunction of functions on $E$
such that the evaluation map at any site $y\in E$ is a continuous
function on the completed space. The following theorem proved by
Aronszajn \cite{Ar} gives a necessary and sufficient condition for
the functional completion.
\begin{thm}[Aronszajn]\label{thm:functional_completion}
Let $\sf F$ be a class of functions on $E$ forming a pre-Hilbert
space. In order that there exists a functional completion of $\sf
F$, it is necessary and sufficient that
\begin{itemize}
\item[(i)] for every fixed $y\in E$, the linear functional $f(y)$ defined in $\sf F$ is continuous; \item[(ii)]
for a Cauchy sequence $\{f_n\}\sbs \sf F$, the condition $f_n(y)\to 0$ for every $y$ implies that $f_n$ itself
converges to $0$ in norm.
\end{itemize}
If the
functional completion is possible, it is unique.
\end{thm}
In our setting, the incomplete class of functions is $\mcH_0$ equipped with the inner product $(\cdot,\cdot)_-$.
We shall demand $\mcH_-$ to be functionally completed. We state all the conditions we need as a hypothesis:

\medskip
(H) The Hermitian positive definite linear operator $A$ on $\mcH_0$ is bounded and satisfies (i)
$\text{ker}A=\{0\}$; (ii) $\mcH_-$ is functionally completed.

\medskip
In the Appendix we will consider some examples of the operators $A$ that satisfy the conditions in (H).

Now $\mcH_-$ being functionally completed, it satisfies, by
definition, that for every $y\in E$, the functional $f(y)$ is
continuous on $\mcH_-$. Notice that by the dual pairing $_-\langle
\cdot,\cdot\rangle_+$, it is equivalent to saying that $e_y\in
\mcH_+$ for any $y\in E$. In fact, it is not hard to check that
the functional $_-\langle  \cdot,g\rangle_+$ on $\mcH_-$ has norm
$\|g\|_+$ for any $g\in \mcH_+$, and the functional $_-\langle
f,\cdot\rangle_+$ on $\mcH_+$ has norm $\|f\|_-$ for each $f\in
\mcH_-$. Moreover, by the isometries $A:\mcH_-\to \mcH_+$ and its
inverse $A^{-1}:\mcH_+\to \mcH_-$, it is easy to check that
\begin{equation}\label{eq:dual_pair_another_form}
_-\langle  \cdot,g\rangle_+=(\cdot,A^{-1}g)_- \text{ and
}_-\langle f,\cdot\rangle_+=(Af,\cdot)_+,\quad f\in \mcH_-,\, g\in
\mcH_+.
\end{equation}
That is, $\mcH_+$ and $\mcH_-$ are respectively the dual spaces of each other via the dual pairing $_-\langle
\cdot,\cdot\rangle_+$. Now if $e_y\in \mcH_+$ for every $y\in E$, then obviously $f(y)=\,_-\langle
f,e_y\rangle_+$ is continuous on $\mcH_-$. On the other hand, suppose that the functional $f(y)$ is continuous
on $\mcH_-$ for every $y\in E$. Then, for each $y$, by the above observation, there is a unique element $l_y\in
\mcH_+$ such that
\begin{equation}\label{eq:linear_functional}
f(y)=\, _-\langle  f,l_y\rangle_+,\quad f\in \mcH_-.
\end{equation}
Since finitely supported vectors $f$ are dense in $\mcH_-$ and for
those vectors $f$ we have $_-\langle
f,l_y\rangle_+=\sum_x\ol{f(x)}l_y(x)$, $l_y$ must be $e_y$.

Finally, we notice that since for any fixed $y\in E$ the
functionals $\mcH_-\ni f\mapsto f(y)$ and $\mcH_+\ni g\mapsto
g(y)$ are continuous, respectively in $\mcH_-$ and $\mcH_+$, it is
obvious that
\begin{equation}\label{eq:dual_pair_realization}
_-\langle  f,g\rangle_+=\sum_{x\in E}\ol{f(x)}g(x), \text{ if
either }f \text{ or }g \text{ is locally supported}.
\end{equation}

\subsection{Determinantal Point Processes on Discrete Sets}
Determinantal point processes, or fermion random point fields, are probability measures on the configuration
space of, say, particles. The particles may move on the continuum spaces or on the discrete spaces. In this
paper we will focus on the DPP's on the discrete sets.

The basics of DPP's including their definitions and basic properties can be found in several papers \cite{GY, L,
LS, M, ST1, ST2, SY, So}. We will review the definition of DPP's mainly from the paper \cite{ST2}. Let $E$ be a
countable set and let $K$ be a Hermitian positive definite bounded linear operator on the Hilbert space
$\mcH_0=l^2(E)$. Let $\mcX$ be the configuration space on $E$, that is, $\mcX$ is the class of all subsets of
$E$. We frequently understand a point $\x=(x_i)_{i=1,2,\cdots}\in \mcX$ as a configuration of particles located
at the sites $x_i\in E$, $i=1,2,\cdots$. The following theorem gives an existence theorem for DPP's. We state it
as appeared in \cite{ST2}.
\begin{thm}\label{thm:existence_of_DPP}
Let $E$ be a countable discrete space and $K$ be a Hermitian bounded operator on $\mcH_0=l^2(E)$. Assume that
$0\le K\le I$. Then, there exists a unique probability Borel measure $\m$ on $\mcX$ such that for any finite
subset $X\sbs E$,
\begin{equation}\label{eq:correlation_function}
\m(\{\x\in \mcX:\,\x\sps X\})=\det(K(x,y))_{x,y\in X}.
\end{equation}
\end{thm}
The $\s$-algebra on $\mcX$ is induced from the product topology on
$\{0,1\}^E$ (see Section 4). Here we remark that the left hand
side of \eqref{eq:correlation_function} is just the correlation
function of the probability measure $\m$, thus the theorem says
that the correlation functions of DPP's are given by the
determinants of positive definite kernel functions.

The most useful feature in the theory of DPP's is that there can
be given an exact formula for the density functions of local
marginals. For each subset $\L\sbs E$, let $P_\L$ denote the
projection operator on $\mcH_0$ onto the space of vectors which
have supports on the set $\L$. Let $K_\L:=P_\L K P_\L$ be the
restriction of $K$ on the projection space. Given a configuration
$\x\in \mcX$, we let $\x_\L$ be the restriction of $\x$ on the set
$\L$, i.e.,
\begin{equation}\label{eq:restriction}
\x_\L:=\x\cap \L.
\end{equation}
For each finite subset $\L\sbs E$, assuming first that $I_\L-K_\L$
is invertible, we define
\begin{equation}\label{eq:local_interaction_operator}
A_{[\L]}:=K_\L(I_\L-K_\L)^{-1}.
\end{equation}
Then for the DPP $\m$ corresponding to the operator $K$, the marginals are given by the formula: for each finite
subset $\L\sbs E$ and fixed $\x\in \mcX$,
\begin{equation}\label{eq:marginals}
\m(\{\z:\z_\L=\x_\L\})=\det(I_\L-K_\L)\det(A_{[\L]}(x,y))_{x,y\in
\x_\L},
\end{equation}
where $A_{[\L]}(x,y)$, $x,y\in \L$, denotes the matrix components of $A_{[\L]}$. Though in this paper we will
confine ourselves to the case where $A_{[\L]}$ is well-defined as a bounded operator, we remark that the formula
\eqref{eq:marginals} is meaningful even if $K_\L$ has $1$ in its spectrum \cite{ST2, So}.

\subsection{Results}
First we will consider a variational principle for the positive
definite operator $A$ introduced in Subsection 2.1. Since we are
assuming that $\mcH_-$ is functionally completed, for any $y\in E$
the functional $f(y)$ is continuous on $\mcH_-$, or $e_y\in
\mcH_+$. This condition, on the other hand, is equivalent to the
one that $\mcH_-$ is a RKHS \cite[p 343]{Ar}. Let $B(x,y)$ be the
RK for $\mcH_-$. From the reproducing property we see that
$B(x,y)$ is the value of the function $A^{-1}e_y$ at $x$ \cite[p
344]{Ar}, that is
\begin{equation}\label{eq:kernel_B_recognization}
B(x,y)=\,_+\langle  e_x,A^{-1}e_y\rangle_-,\quad x,y\in E.
\end{equation}
Let $x_0\in E$ be a fixed point and let $R_1$ and $R_2$ be any two
subsets of $E$ such that $E$ is partitioned into three sets:
\begin{equation}\label{eq:partition_of_E}
E=\{x_0\}\cup R_1\cup R_2.
\end{equation}
For each $\D\sbs E$, we let ${\sf F}_{\text{loc},\D}$ be the local functions supported on $\D$:
\begin{equation}\label{eq:local_functions}
{\sf F}_{\text{loc},\D}:=\text{the class of finite linear
combinations of }\{e_x:\,x\in \D\}.
\end{equation}
In the sequel, we denote by $\L\Subset E$ that $\L$ is a finite
subset of $E$. We are concerned with the following numbers. For
each $\L\ssbs E$, define
\begin{equation}\label{eq:alpha_local} \a_\L:=\inf_{f\in
{\sf F}_{\text{loc},\L\cap R_1}}\|e_{x_0}-f\|_-^2
\end{equation}
and
\begin{equation}\label{eq:beta_local}
\b_\L:=\inf_{g\in {\sf F}_{\text{loc},\L\cap R_2}}\|e_{x_0}-g\|_+^2.
\end{equation}
Obviously, both of the sequences of nonnegative numbers
$\{\a_\L\}_{\L\ssbs E}$ and $\{\b_\L\}_{\L\ssbs E}$ decrease as
$\L$ increases. We let
\begin{equation}\label{eq:alpha_and_beta}
\a:=\lim_{\L\ua E}\a_\L\quad \text{and}\quad \b:=\lim_{\L\ua
E}\b_\L.
\end{equation}
One of the main result of this paper is the following:
\begin{thm}\label{thm:variational_principle}
Let the operator $A$ satisfy the conditions in the hypothesis (H).
Then the product of the numbers $\a$ and $\b$ defined in
\eqref{eq:alpha_and_beta} is one: $\a\b=1$.
\end{thm}
We remark that the result of the theorem was obtained by Shirai
and Takahashi \cite[Theoem 6.3]{ST2} in the case that the bounded
operator $A$ is strictly greater than $0$, i.e., $0<cI\le A$ for
some positive constant $c$.

One of the main purpose of this paper is to apply the above result
to show the Gibbsianness of some DPP's. Let $A$ be an operator on
$\mcH_0$ that satisfies the hypothesis (H). Let $\m$ be the DPP
corresponding to the operator $K:=A(I+A)^{-1}$. Given a fixed
point $x_0\in E$ and a configuration $\x\in \mcX$ with $x_0\notin
\x$, and for each $\L\ssbs E$, let $\a_{[\L]}$ be the conditional
probability of finding a particle at the site $x_0$ given the
particle configuration $\x_\L$ in $\L$:
\begin{equation}\label{eq:conditional_probability}
\a_{[\L]}:=\m(x_0\x_\L|\x_\L)=\frac{\m_\L(x_0\x_\L)}{\m_\L(\x_\L)},
\end{equation}
where we have simplified the event $\{\z\in
\mcX:\,\z_\L=\x_\L\}\equiv \x_\L$,  etc, and $x_0\x_\L=\{x_0\}\cup
\x_\L$. By \eqref{eq:marginals}, $\a_{[\L]}$ is computed via the
ratio of determinants:
\begin{equation}\label{eq:alpha_local_by_determinants}
\a_{[\L]}=\frac{\det A_{[\L]}(x_0\x_\L,x_0\x_\L)}{\det
A_{[\L]}(\x_\L,\x_\L)},
\end{equation}
where $A_{[\L]}=K_\L(I_\L-K_\L)^{-1}$ and $A_{[\L]}
(\x_\L,\x_\L)=(A_{[\L]}(x,y))_{x,y\in \x_\L}$. We are interested
in the behavior of the sequence $\{\a_{[\L]}\}$ as $\L$ increases
to $E$. The following theorem gives the answer.
\begin{thm}\label{thm:global_PI}
Let the operator $A$ satisfy the conditions in (H). Then
\begin{equation}\label{eq:gloval_PI}
\lim_{\L\ua E}\a_{[\L]}=\a,
\end{equation}
where $\a$ is given in \eqref{eq:alpha_local} and \eqref{eq:alpha_and_beta} with $R_1=\x$ and $R_2=E\sm(\x\cup
\{x_0\})$.
\end{thm}
A corollary to this theorem is that the DPP $\m$ corresponding to
the operator  $A(I+A)^{-1}$ is a Gibbs measure. We state this as a
theorem.
\begin{thm}\label{thm:Gibbs_measure}
Let the operator $A$ satisfy the conditions in (H). Then the DPP $\m$  corresponding to the operator
$A(I+A)^{-1}$ is a Gibbs measure. The interaction potential is given by the logarithm of determinants of
submatrices of $A$: for any finite configuration $\x\in \mcX$, the interaction potential $V(\x)$ is
\begin{equation}\label{eq:potential_energy}
V(\x)=-\log\det(A(x,y))_{x,y\in \x}.
\end{equation}
Moreover, $\m$ is the only Gibbs measure for the potential energy \eqref{eq:potential_energy}.
\end{thm}
The above result also extends that obtained in \cite[Theorem
6.2]{ST2},  where $K\equiv A(I+A)^{-1}$ is assumed to have its
spectrum in the open interval $(0,1)$. We also notice that the
idea developed in refs. \cite{GY} and \cite{Y}, which concerns
exclusively with continuum models, can be applied to discrete
model and would get some result on the Gibbsianness of $\m$. The
result would look like the following (cf. \cite[Proposition
3.9]{GY}): Let $E\equiv \mbZ^d$ and suppose that (i) $A$ is of
finite range in the sense that $A(x,y)=0$ if $|x-y|\ge R$ for some
finite number $R>0$ and (ii) $\m$ does not percolate. Then $\m$ is
a Gibbs measure corresponding to the potential in
\eqref{eq:potential_energy}. Our result \ref{thm:Gibbs_measure} is
stronger than this, too.

\section{Proof of the Variational Principle}
In this Section we prove Theorem \ref{thm:variational_principle}.
The most important tool in the proof is the theory of restrictions
and projections in the RKHS's. In Subsection 3.1, we deal with the
variational principle in the finite systems. In Subsection 3.2, we
first introduce the restriction theory in the RKHS's and then
discuss the limit theorems of RK's. The proof of Theorem
\ref{thm:variational_principle} is given in Subsection 3.3.

\subsection{Variational Principle in the Finite Systems}

We discuss the variational principle for positive definite matrices on a finite set. Let $\L\ssbs E$ be a finite
set and let $(C(x,y))_{x,y\in \L}$ be a positive definite matrix with an inverse $C^{-1}$. We define two norms
on the class ${\sf F}_\L$ of functions on $\L$ as follows:
\begin{equation}\label{eq:-norm_on_local}
\|f\|_-^2:=\sum_{x,y\in \L}\ol{f(x)}C(x,y)f(y),\quad f\in {\sf F}_\L
\end{equation}
and
\begin{equation}\label{eq:+norm_on_local}
\|g\|_+^2:=\sum_{x,y\in \L}\ol{g(x)}C^{-1}(x,y)g(y),\quad g\in {\sf F}_\L.
\end{equation}
Suppose that $\L=\{x_0\}\cup\L_1\cup\L_2$ is a partition of $\L$
with disjoint sets $\{x_0\}$, $\L_1$, and $\L_2$. Similarly to
\eqref{eq:alpha_local}-\eqref{eq:beta_local} we define
\begin{equation}\label{eq:a_local}
a:=\inf_{f\in {\sf F}_{\L_1}}\|e_{x_0}-f\|_-^2
\end{equation}
and
\begin{equation}\label{eq:b_local}
b:=\inf_{g\in {\sf F}_{\L_2}}\|e_{x_0}-g\|_+^2.
\end{equation}
In the above ${\sf F}_{\L_i}$ denotes the class of functions on
$\L_i$, $i=1,\,2$. Applying the method of finding extreme values
of functions of several variables and using the elementary
properties of determinants of finite matrices, we obtain the
following results, which, as a matter of fact, take a role of
recipe for the theory in the infinite systems (cf. \cite[Section
6]{ST2}). Below we denote by $C(\L_1,\L_2)$ the submatrix
$(C(x,y))_{x\in \L_1,\,y\in \L_2}$ for any subsets $\L_1$ and
$\L_2$ of $\L$. We also simplify $\{x\}\cup \L_1$ by $x\L_1$ for
$x\notin \L_1$.
\begin{prop}\label{prop:results_of_local_systems}
Let $(C(x,y))_{x,y\in \L}$ be a Hermitian positive definite matrix
on a finite set $\L$ with inverse $C^{-1}$. Let
$\L=\{x_0\}\cup\L_1\cup\L_2$ be a partition of $\L$ and let the
norms $\|\cdot\|_-$ and $\|\cdot\|_+$, and the numbers $a$ and $b$
be defined as in \eqref{eq:-norm_on_local}-\eqref{eq:b_local}.
Then the following results hold:
\begin{itemize} \item[(a)] The minimum values $a$ and $b$ are
attained respectively at the unique vectors
$f_0=C(\L_1,\L_1)^{-1}C(\L_1,x_0)$ and
$g_0=(C^{-1}(\L_2,\L_2))^{-1}C^{-1}(\L_2,x_0)$:
\begin{equation}\label{eq:a_and_b}
a=\|e_{x_0}-f_0\|_-^2;\quad b=\|e_{x_0}-g_0\|_+^2.
\end{equation}
\item[(b)]
\begin{eqnarray}\label{eq:a_2}
a&=&\frac{\det C(x_0\L_1,x_0\L_1)}{\det
C(\L_1,\L_1)}=(C(x_0\L_1,x_0\L_1)^{-1}(x_0,x_0))^{-1}\nonumber\\
&=& C(x_0,x_0)-C(x_0,\L_1)C(\L_1,\L_1)^{-1}C(\L_1,x_0)
\end{eqnarray}
and similarly
\begin{eqnarray}\label{eq:b_2}
b&=&\frac{\det C^{-1}(x_0\L_2,x_0\L_2)}{\det
C^{-1}(\L_2,\L_2)}=((C^{-1}(x_0\L_2,x_0\L_2))^{-1}(x_0,x_0))^{-1}\nonumber\\
&=&
C^{-1}(x_0,x_0)-C^{-1}(x_0,\L_2)(C^{-1}(\L_2,\L_2))^{-1}C^{-1}(\L_2,x_0)
\end{eqnarray}
\item[(c)] $ab=1$.
\end{itemize}
\end{prop}

\subsection{Restrictions in RKHS's and Limit Theorems of RK's}

In this Subsection, we discuss the restriction and projection
theories in RKHS's and the limit theorems of RK's. These are
crucial to characterize the values $\a$ and $\b$ in
\eqref{eq:alpha_and_beta} more concretely. The results we need
have been already obtained in \cite{Ar}. For the readers'
convenience, however, we provide it here.

Let us begin with an introduction of the restriction theory in the
RKHS's. Suppose that $\mcH$ is a RKHS with kernel $K(x,y)$,
$x,y\in E$. $\mcH$ might be $\mcH_-$ or $\mcH_+$ of our concern.
For each subset $\L\sbs E$, the function $K_\L(x,y)$, the
restriction of $K(x,y)$ to $\L$, is still positive definite. The
following theorem was proved by Aronszajn \cite[p 351]{Ar}:
\begin{thm}\label{thm:Aronszajn}
The function $K(x,y)$ restricted to a subset $\L\sbs E$ is the reproducing kernel of the class $\mcH_\L$ of all
restrictions of functions of $\mcH$ to the subset $\L$. For any such restriction $f_\L\in \mcH_\L$, the norm
$\|f_\L\|_\L$ is the minimum of $\|f\|$ (the norm of $f$ in $\mcH$) for all $f\in \mcH$ whose restrictions to
$\L$ are $f_\L$.
\end{thm}
When it is needed to designate the kernel, we use the notations
$\mcH_{\L;K}$ and $\|\cdot\|_{\L;K}$ respectively for the
restriction spaces and norms. The basic argument in Theorem
\ref{thm:Aronszajn} is the following. First let ${\sf F}^0\sbs
\mcH$ be the class of functions that vanish on $\L$. This is a
closed subspace and let ${\sf F}':=\mcH\ominus{\sf F}^0$ be the
orthogonal complement of ${\sf F}^0$. It is not hard to show that
all the functions $f\in \mcH$ which have the same restriction
$f_\L$ on $\L$ have a common projection $f'$ on ${\sf F}'$ and
that the restriction of $f'$ to $\L$ is equal to $f_\L$. Clearly,
among all these functions $f$, $f'$ is the one which has the
smallest norm. We define
\begin{equation}\label{eq:norm_on_restriction}
\|f_\L\|_\L:=\|f'\|.
\end{equation}
The norm $\|\cdot\|_\L$ on $\mcH_\L$ defined this way is the one stated in the theorem. We refer to \cite[p
351]{Ar} for the details.

Next we discuss the limit theorems of RK's. We will consider two
kinds of limits.

{\it A. The case of decreasing sequence.} Let $\{E_n\}$ be an
increasing sequence of sets with $E=\cup_{n=1}^\infty E_n$. For
each $n=1,2,\cdots$, let ${\sf F}_n$ be a RKHS defined in $E_n$
with RK $K_n(x,y)$, $x,y\in E_n$. we denote the norm in the space
${\sf F}_n$ by $\|\cdot\|_n$, $n\ge 1$. For a function $f_n\in
{\sf F}_n$ we will denote by $f_{nm}$, $m\le n$, the restriction
of $f_n$ to the set $E_m\sbs E_n$. We shall suppose the following
two conditions:
\begin{enumerate}\item[(A1)] for every $f_n\in {\sf F}_n$ and every $m\le n$, $f_{nm}\in {\sf F}_m$; \item[(A2)] for
every $f_n\in {\sf F}_n$ and every $m\le n$, $\|f_{nm}\|_m\le
\|f_n\|_n$.
\end{enumerate}
From (A2) we see by \cite[Theorem II of Section 7]{Ar} that
\begin{equation}\label{eq:monotone_decreasing}
K_{nm}\ll K_m,\quad m<n,
\end{equation}
meaning that $K_m(x,y)-K_{nm}(x,y)$, $x,y\in E_m$, is a positive
definite function, where $K_{nm}$ is the restriction of $K_n$ to
the set $E_m$. The following theorem appears in \cite[Theorem I,
Section 9]{Ar}:
\begin{thm}\label{thm:AronszajnI} Under the above assumptions on the classes ${\sf F}_n$, the kernels
$K_n$ converge to a Kernel $K_0(x,y)$defined for all $x,y$ in $E$.
$K_0$ is the RK of the class ${\sf F}_0$ of all functions $f_0$
defined in $E$ such that \begin{enumerate}\item[(i)] their
restrictions $f_{0n}$ in $E_n$ belong to ${\sf F}_n$,
$n=1,2,\cdots$; \item[(ii)] $\lim_{n\to
\infty}\|f_{0n}\|_n<\infty$.
\end{enumerate} The norm of $f_0\in {\sf F}_0$ is given by $\|f_0\|_0=\lim_{n\to \infty}\|f_{0n}\|_n$.
\end{thm}

{\it B. The case of increasing sequence.} Let $\{E_n\}$ be a
decreasing sequence of sets and $R$ be their intersection:
\begin{equation}
R=\cap_{n=1}^\infty E_n.
\end{equation}
As in the case A, let ${\sf F}_n$, $n=1,2\cdots$, be the RKHS's
with corresponding kernel functions $K_n(x,y)$, $x,y\in E_n$,
$n\ge 1$. As before, we define the restrictions $f_{nm}$ for
$f_n\in {\sf F}_n$, but now $m$ has to be greater than $n$. We
suppose that ${\sf F}_n$ form an increasing sequence and the norms
$\|\cdot\|_n$ form a decreasing sequence satisfying the following
two conditions:
\begin{enumerate}\item[(B1)] for every $f_n\in {\sf F}_n$ and every
$m\ge n$, $f_{nm}\in {\sf F}_m$; \item[(B2)] for every $f_n\in
{\sf F}_n$ and every $m\ge n$, $\|f_{nm}\|_m\le \|f_n\|_n$.
\end{enumerate}
We then get for the restrictions $K_{nm}$ of $K_n$ the formula
\begin{equation}
K_{nm}\ll K_m,\quad \text{for }m\ge n.
\end{equation}
For each $y\in R$, $\{K_m(y,y)\}$ is an increasing sequence of positive numbers. Its limit may be infinite. We
define, consequently,
\begin{equation}
R_0:= \text{ the set of }y\in R \text{ such that
}K_0(y,y):=\lim_{m\to \infty} K_m(y,y)<\infty.
\end{equation}
Suppose that $R_0$ is not empty and let ${\sf F}_0$ be the class
of all restrictions $f_{n0}$ of functions $f_n\in {\sf F}_n$
($n=1,2,\cdots$) to the set $R_0$. From (B2), the limit
$\lim_{k\to \infty}\|f_{nk}\|_k$ exists and we define a norm
$\|\cdot\|_0^\sim$ on ${\sf F}_0$ by\footnote{The original
definition in \cite{Ar} is such that $\|f_{n0}\|_0:=\lim_{k\to
\infty}\|f_{nk}\|_k$, but it seems that there is no way to
guarantee that $\|f_{n0}\|_0=\|g_{n0}\|_0$ for different $f_n$ and
$g_n$ in ${\sf F}_n$ with $f_{n0}=g_{n0}$. However, all the
arguments in \cite{Ar} hold true even if the new norm
$\|\cdot\|_0^\sim$ in \eqref{eq:0*_norm} is used. In particular,
the Theorem \ref{thm:AronszajnII} below holds.}
\begin{equation}\label{eq:0*_norm}
\|f\|_0^\sim:=\inf\lim_{k\to \infty}\|f_{nk}\|_k,\quad f\in {\sf
F}_0,
\end{equation}
where the infimum is taken over all functions $f_n\in {\sf F}_n$,
$n\ge 1$, whose restrictions to $R_0$ are $f$, i.e.,
$f(y)=f_{n0}$, $y\in R_0$, for some $f_n\in {\sf F}_n$. Now we
construct a new space ${\sf F}_0^*$ and norm $\|\cdot\|_0^*$ on
it. Let ${\sf F}_0^*$ be the class of all functions $f_0^*$ on
$R_0$ such that there is a Cauchy sequence $\{f_0^{(n)}\}\sbs {\sf
F}_0$ satisfying
\begin{equation}\label{eq:vectors_of_F0*}
f_0^*(x)=\lim_{n\to \infty}f_0^{(n)}(x),\quad \text{for all }x\in
R_0.
\end{equation}
For those vectors $f_0^*$ we define a norm
\begin{equation}
\|f_0^*\|_0^*:=\min\lim_{n\to \infty}\|f_0^{(n)}\|_0^\sim,
\end{equation}
the minimum being taken over all Cauchy sequences
$\{f_0^{(n)}\}\sbs {\sf F}_0$ satisfying
\eqref{eq:vectors_of_F0*}. There exists at least one Cauchy
sequence for which the minimum is attained. Such sequences are
called determining $f_0^*$. The scalar product corresponding to
$\|\cdot\|_0^*$ is defined by
\begin{equation}
(f_0^*,g_0^*)_0^*:=\lim_{n\to \infty}(f_0^{(n)},g_0^{(n)})_0^\sim
\end{equation}
for any two Cauchy sequences $\{f_0^{(n)}\}$ and $\{g_0^{(n)}\}$
determining $f_0^*$ and $g_0^*$, respectively. We refer to
\cite[Section 9]{Ar} for the details. The following theorem is in
\cite[Theorem II, Section 9]{Ar}:
\begin{thm}\label{thm:AronszajnII}
In the setting of the case B, the restrictions $K_{n0}(x,y)$ for
every fixed $y\in R_0$ form a Cauchy sequence in ${\sf F}_0$. They
converge to a function $K_0^*(x,y)\in {\sf F}_0^*$ which is the RK
of ${\sf F}_0^*$.
\end{thm}

As an application of Theorem \ref{thm:AronszajnII}, we prove the convergence of norms in the perturbed RKHS's,
which will be used in the proof of Theorem \ref{thm:variational_principle}. Let $A$ be the operator of our
concern satisfying the conditions in the hypothesis (H). For each $\e>0$ we define new operators as follows:
\begin{equation}\label{eq:perturbed_operators}
A(\e):=A+\e\quad \text{and}\quad B(\e):=A(\e)^{-1},\quad \e>0.
\end{equation}
Let $R\sbs E$ be any subset of $E$. Following Theorem
\ref{thm:Aronszajn}, we let $\|\cdot\|_{R;B}$ be the norm of the
RKHS $\mcH_{R;B}$ consisting of all restrictions of vectors in
$\mcH_-$ to the set $R$ and having a RK $B_R(x,y)$, $x,y\in R$,
the restriction of $B(x,y)$ to the set $R$. Similarly,
$\|\cdot\|_{R;B(\e)}$ denotes the norm defined by replacing $B$
with $B(\e)$. We want to prove the convergence $\|f\|_{R;B(\e)}\to
\|f\|_{R;B}$ for all $f\in l^2(R)$ as $\e\to 0$. See Lemma
\ref{lem:limit_of_perturbation}. For that purpose we proceed as
follows. Let $\mcH_{R;B}'\sbs l^2(R)$ be the dual space of
$\mcH_{R;B}$: an element $g\in l^2(R)$ belongs to $\mcH_{R;B}'$ if
and only if the (anti-)linear functional
\begin{equation}\label{eq:anti_linear_functional}
{}_{R;B}\langle  f,g\rangle_{R;B}':=\sum_{x\in
R}\ol{f(x)}g(x),\quad f\in l^2(R),
\end{equation}
is continuous w.r.t. $\|\cdot\|_{R;B}$-norm, i.e., there exists $M(g)>0$ such that
\begin{equation}
|{}_{R;B}\langle  f,g\rangle_{R;B}'|\le M(g)\|f\|_{R;B},\quad
\text{for all }f\in l^2(R).
\end{equation}
For each $g\in \mcH_{R;B}'$ we extend the functional of
\eqref{eq:anti_linear_functional} to the whole space
$\mcH_{R;B}\sps l^2(R)$ and keep the dual pairing notation
${}_{R;B}\langle  \cdot,\cdot\rangle_{R;B}'$. We denote the norm
in $\mcH_{R;B}'$ by $\|\cdot\|_{R;B}'$. As in the case of the dual
pairing ${}_-\langle  \cdot,\cdot\rangle_+$ we see that for any
$f\in \mcH_{R;B}$, $(B_R)^{-1}f\in \mcH_{R;B}'$ and
\begin{equation}
\|(B_R)^{-1}f\|_{R;B}'=\|f\|_{R;B}.
\end{equation}
It is not hard to show that for any $h\in l^2(R)$,
\begin{equation}\label{eq:norm_order}
(B_R)^{-1}h\in \mcH_+ \text{ and }\|(B_R)^{-1}h\|_+\le \|(B_R)^{-1}h\|_{R;B}'.
\end{equation}
In fact, we have for any $f\in \mcH_0=l^2(E)$,
\begin{eqnarray}\label{eq:proof_of_anti_linear_functional}
|{}_-\langle  f,(B_R)^{-1}h\rangle_+|&=&|\sum_{x\in R}\ol{f(x)}(B_R)^{-1}h(x)|\nonumber\\
&=&|{}_{R;B}\langle  f_R,(B_R)^{-1}h\rangle_{R;B}'|\nonumber\\
&\le &\|f_R\|_{R;B}\|(B_R)^{-1}h\|_{R;B}'\nonumber\\
&\le &\|f\|_-\|(B_R)^{-1}h\|_{R;B}',
\end{eqnarray}
where $f_R$ is the restriction of $f$ to $R$. Since $\mcH_0$ is dense in $\mcH_-$,
\eqref{eq:proof_of_anti_linear_functional} proves \eqref{eq:norm_order}. Because $l^2(R)$ is dense in
$\mcH_{R;B}$, \eqref{eq:norm_order} also shows that
\begin{equation}
\mcH_{R;B}'\sbs\mcH_+\cap l^2(R).
\end{equation}
\begin{lem}\label{lem:vanishing_for_disjoint_vectors}
Let $g\in \mcH_{R;B}'\sbs \mcH_+\cap l^2(R)$. Then for any $f\in
\mcH_-$ that vanishes on $R$, we have ${}_-\langle
f,g\rangle_+=0$.
\end{lem}
\medskip
\Proof Denote by $P_R$ the restriction operator $P_R:\mcH_-\to \mcH_{R;B}$ defined by $P_Rf:=f_R$ for all $f\in
\mcH_-$. Since $\|f_R\|_{R;B}\le \|f\|_-$, the operator $P_R$ is bounded with norm less than or equal to $1$.
Now let $g$ and $f$ be as in the statement of the lemma. Let $\{f_n\}$ be any sequence in $\mcH_0$ that
converges to $f$ in $\mcH_-$. Then since $g\in l^2(R)$, by using the continuity of the operator $P_R$ in
$\mcH_-$, we have
\begin{eqnarray*}
{}_-\langle  f,g\rangle_+&=&\lim_{n\to \infty}{}_-\langle  f_n,g\rangle_+\\
&=&\lim_{n\to \infty}\sum_{x\in R}\ol{f_n(x)}g(x)\\
&=&\lim_{n\to \infty}{}_{R;B}\langle  P_Rf_n,g\rangle_{R;B}'\\
&=& {}_{R;B}\langle  P_Rf,g\rangle_{R;B}'\\
&=&0,
\end{eqnarray*}
because $P_Rf=0$. \EndProof
\begin{lem}\label{lem:equivalence_of_norm}
For any $h\in l^2(R)$, $\|(B_R)^{-1}h\|_+=\|(B_R)^{-1}h\|_{R;B}'$.
\end{lem}
\medskip
\Proof By \eqref{eq:norm_order} it is enough to show that $\|(B_R)^{-1}h\|_{R;B}'\le\|(B_R)^{-1}h\|_+$. We have
\begin{eqnarray}\label{eq:another_form}
(\|(B_R)^{-1}h\|_{R;B}')^2&=&\|h\|_{R;B}^2\nonumber\\
&=&{}_{R;B}\langle  h,(B_R)^{-1}h\rangle_{R;B}'\nonumber\\
&=&{}_-\langle  h,(B_R)^{-1}h\rangle_+.
\end{eqnarray}
Let $h'\in \mcH_-$ be the element such that $P_Rh'=h$ and $\|h'\|_-=\|h\|_{R;B}$ (see
\eqref{eq:norm_on_restriction}). Notice that $h'-h\in \mcH_-$ vanishes on $R$ and $(B_R)^{-1}h\in \mcH_{R;B}'$.
Thus by \eqref{eq:another_form} and Lemma \ref{lem:vanishing_for_disjoint_vectors} we have
\begin{eqnarray*}
\|h\|_{R;B}^2&=&{}_-\langle  h,(B_R)^{-1}h\rangle_+\\
&=&{}_-\langle  h',(B_R)^{-1}h\rangle_+\\
&\le &\|h'\|_-\|(B_R)^{-1}h\|_+\\
&=&\|h\|_{R;B}\|(B_R)^{-1}h\|_+.
\end{eqnarray*}
This, together with \eqref{eq:another_form}, proves that $\|(B_R)^{-1}h\|_{R;B}'\le \|(B_R)^{-1}h\|_+$.
\EndProof

Recall the definition $A(\e)=A+\e$ and $B(\e)=A(\e)^{-1}$ for
$\e>0$.
\begin{lem}\label{lem:limit_of_perturbation}
Let $R\sbs E$ be any set. Then for any $f\in l^2(R)$,
\begin{equation}\label{eq:limit_of_purturbation}
\lim_{\e\to 0}\|f\|_{R;B(\e)}=\|f\|_{R;B}.
\end{equation}
\end{lem}
\medskip
\Proof First we show that
\begin{equation}\label{eq:pointwise_convergence}
\lim_{\e\to 0}B(\e)(x,y)=B(x,y),\quad\text{for all }x,y\in E.
\end{equation}
It is obvious that
\begin{equation}
B(\e)\ll B(\e')\quad \text{for }0<\e'<\e,
\end{equation}
in the sense defined in \eqref{eq:monotone_decreasing}. Also, it holds trivially that
\begin{equation}
B(\e)(y,y)\le B(y,y)<\infty,\quad \forall\,\,y\in E.
\end{equation}
Moreover, for each fixed $\e>0$, since $B(\e)$ is bounded and strictly positive, the norms $\|\cdot\|_{-;\e}$
($:=\|\cdot\|_{E;B(\e)}$) and $\|\cdot\|_0$ are equivalent on $\mcH_0=l^2(E)$. That is, as a set, $\mcH_{-;\e}$
($:=\mcH_{E;B(\e)}$) is the same as $\mcH_0$.

 Now for each $f\in \mcH_0$, the norm $\|f\|_{-;\e}$ decreases as $\e$ decreases. It is easy to check
 that
 \begin{equation}
 \lim_{\e\to 0}\|f\|_{-;\e}=\|f\|_-.
 \end{equation}
 In fact, for $f\in \mcH_0$,
 \begin{eqnarray*}
 \lim_{\e\to 0}\|f\|_{-;\e}^2&=&\lim_{\e\to 0}(f,A(\e)f)_0\\
 &=&(f,Af)_0\\
 &=&\|f\|_-^2.
 \end{eqnarray*}
Since the norm $\|\cdot\|_-$ is the one for the RKHS with kernel $B(x,y)$, the equality
\eqref{eq:pointwise_convergence} follows from Theorem \ref{thm:AronszajnII} (see the remark on \cite[p
368]{Ar}).

Let us now prove \eqref{eq:limit_of_purturbation}. Obviously, for each $f\in l^2(R)$, $\|f\|_{R;B(\e)}$
decreases as $\e$ decreases and $\|f\|_{R;B(\e)}\ge \|f\|_{R;B}$ for all $\e>0$. Thus the limit
\begin{equation}
\|f\|_0^\sim :=\lim_{\e\to 0}\|f\|_{R;B(\e)},\quad f\in l^2(R),
\end{equation}
defines a norm on $l^2(R)$. Now we have to show
$\|f\|_0^\sim=\|f\|_{R;B}$. Considering the dual norms it is
equivalent to showing that
\begin{equation}\label{eq:equivalent_condition}
\lim_{\e\to 0}\|g\|_{R;B(\e)}'=\|g\|_{R;B}'
\end{equation}
for $g\in l^2(R)$ whenever the limit is finite. Thus suppose that $g\in l^2(R)$ and $\lim_{\e\to
0}\|g\|_{R;B(\e)}'$ is finite. Since $B(\e)$ is a strictly positive and bounded operator, we see that
\begin{eqnarray}\label{eq:form_representation}
(\|g\|_{R;B(\e)}')^2&=&(g,B(\e)_Rg)_0\nonumber\\
&=&(g,\frac{B}{I+\e B}g)_0.
\end{eqnarray}
Now consider the form $\mcE$ on $\mcH_0$ generated by the operator $B$:
\begin{equation}
\mcE(f,g):=(f,Bg)_0 , \quad f,g\in \text{ran}A.
\end{equation}
Then the space $\mcH_+$ is nothing but the closure of $\text{ran}A$ w.r.t. this form norm. We denote the closure
of the form $(\mcE,\text{ran}A)$ by $(\mcE,D(\mcE))$. Using this notation, the last quantity in
\eqref{eq:form_representation} becomes $\mcE(\frac1{\sqrt{I+\e B}}g,\frac1{\sqrt{I+\e B}}g)$. By the assumption,
these values are bounded from above (as $\e$ varies). On the other hand, $\frac1{\sqrt{I+\e B}}g$ converges
strongly to $g$ as $\e$ goes to $0$. By \cite[Lemma 2.12]{MR}, we conclude that $g\in D(\mcE)$, i.e., $g\in
\mcH_+$, and
\begin{eqnarray}\label{eq:form_bound}
\|g\|_+^2=\mcE(g,g)&\le& \liminf_{\e\to 0}\mcE(\frac1{\sqrt{I+\e B}}g,\frac1{\sqrt{I+\e B}}g)\nonumber\\
&=&\lim_{\e\to 0}(\|g\|_{R;B(\e)})^2.
\end{eqnarray}
To finish the proof we notice that the space $(B_R)^{-1}(l^2(R))$ is dense in $\mcH_{R;B}'$ because $l^2(R)$ is
dense in $\mcH_{R;B}$ and the map $(B_R)^{-1}:\,(l^2(R),\|\cdot\|_{R;B})\to \mcH_{R;B}'$ is an isometry.
Therefore, it is enough to check \eqref{eq:equivalent_condition} for those vectors $g$ of the form
$g=(B_R)^{-1}h$ for some $h\in l^2(R)$. We use Lemma \ref{lem:equivalence_of_norm} and \eqref{eq:form_bound}
inserting $(B_R)^{-1}h$ for $g$. Then we get
\begin{equation}\label{eq:dual_norm_bound}
\|(B_R)^{-1}h\|_{R;B}'=\|(B_R)^{-1}h\|_+\le \lim_{\e\to 0}\|(B_R)^{-1}h\|_{R;B(\e)}'\le \|(B_R)^{-1}h\|_{R;B}'.
\end{equation}
The last inequality in the above comes from the fact that $\|g\|_{R;B(\e)}'\le \|g\|_{R;B}'$ for all $g\in
\mcH_{R;B}'$. Eq. \eqref{eq:dual_norm_bound} says that all the quantities there are equal to each other, and we
have proven \eqref{eq:equivalent_condition}. \EndProof

\subsection{Proof of Theorem \ref{thm:variational_principle}}
We are now ready to prove Theorem \ref{thm:variational_principle}. It will be done in several steps.

\medskip
\Proof[ of Theorem \ref{thm:variational_principle}] {\it Step 1: General facts.} Recall the notation ${\sf
F}_{\text{loc},\L}$, the class of local functions supported on $\L$ for the subsets $\L\sbs E$. Since these
spaces are closed in $\mcH_-$, for each $\L\ssbs E$ there exists a unique element $f_{\L_1,0}\in {\sf
F}_{\text{loc},\L_1}$ such that
\begin{equation}
\a_\L=\inf_{f\in {\sf F}_{\text{loc},\L_1}}\|e_{x_0}-f\|_-^2=\|e_{x_0}-f_{\L_1,0}\|_-^2,
\end{equation}
where $\L_1:= \L\cap R_1$. Let $\mcH_{1,-}$ be the closure (in $\mcH_-$) of $\cup_{\L\ssbs E}{\sf
F}_{\text{loc},\L_1}$. Notice that any vector $f\in \mcH_{1,-}$ vanishes on $R_1^c$, that is, it is supported on
$R_1$. In fact, since for each $\L\ssbs E$, ${\sf F}_{\text{loc},\L_1}\sbs {\sf F}_{R_1^c}^0$, the space of
functions that vanish on $R_1^c$, and ${\sf F}_{R_1^c}^0$ is closed, we have $\mcH_{1,-}=\ol{\cup_{\L\ssbs
E}{\sf F}_{\text{loc},\L_1}}\sbs {\sf F}_{R_1^c}^0$. We also notice that for each $\L\ssbs E$, $f_{\L_1,0}$ is
the projection of $e_{x_0}$ onto the space ${\sf F}_{\text{loc},\L_1}={\sf F}_{\L_1^c}^0$, and as $\L$
increases, $f_{\L_1,0}$ converges to the projection of $e_{x_0}$ onto the space $\mcH_{1,-}$, we call it
$f_{R_1,0}$:
\begin{equation}\label{eq:limit_minimizer}
\lim_{\L\ua E}f_{\L_1,0}=f_{R_1,0}\quad (\text{in }\mcH_-).
\end{equation}
Let us now apply the (extended) operator $A$ to the vector
$e_{x_0}-f_{R_1,0}$. We claim that
\begin{equation}\label{eq:A_operates}
A(e_{x_0}-f_{R_1,0})=\a e_{x_0}+a_2\in \mcH_+,
\end{equation}
where the vector $a_2\in \mcH_+$ is supported on $R_2$. In fact, let $A(e_{x_0}-f_{R_1,0})=a_0 e_{x_0}+a_2\in
\mcH_+$ with $a_2$ being supported on $E\sm\{x_0\}$. Since $f_{R_1,0}$ is the projection of $e_{x_0}$ onto the
space $\mcH_{1,-}$, we have
\begin{equation}\label{eq:orthogonality}
(e_{x_0}-f_{R_1,0},f)_-=0\quad \text{for all }f\in {\sf F}_{\text{loc},R_1}.
\end{equation}
Thus, we have for all $f\in {\sf F}_{\text{loc},R_1}$,
\begin{eqnarray}\label{eq:a_2_vanishing}
0&=&(e_{x_0}-f_{R_1,0},f)_-\nonumber\\
&=&_+\langle  A(e_{x_0}-f_{R_1,0}),f\rangle_-\nonumber\\
&=&_+\langle  a_0 e_{x_0}+a_2,f)\rangle_-\nonumber\\
&=&\sum_{x\in R_1}\ol{a_2(x)}f(x),
\end{eqnarray}
because $f$ is a local function supported on $R_1$ (see \eqref{eq:dual_pair_realization}). Since $f\in {\sf
F}_{\text{loc},R_1}$ is arbitrary, the equation \eqref{eq:a_2_vanishing} proves that $a_2$ vanishes on $R_1$.
Similarly, it is easily checked that
\begin{eqnarray}\label{eq:a_0_realization}
\a&=&\|e_{x_0}-f_{R_1,0}\|_-^2\nonumber\\
&=&_-\langle  e_{x_0}-f_{R_1,0},A(e_{x_0}-f_{R_1,0})\rangle_+\nonumber\\
&=&\lim_{\L\ua E}\,_-\langle  e_{x_0}-f_{\L_1,0},a_0 e_{x_0}+a_2\rangle_+\nonumber\\
&=&\lim_{\L\ua E}a_0=a_0.
\end{eqnarray}
We have shown \eqref{eq:A_operates}. Let us now interchange the roles of $A$, $\|\cdot\|_-$, $R_1$, and $\L_1$
by $A^{-1}$, $\|\cdot\|_+$, $R_2$, and $\L_2$, respectively. Then we have for each $\L\ssbs E$,
\begin{equation}
\b_\L=\inf_{g\in {\sf F}_{\text{loc},\L_2}}\|e_{x_0}-g\|_+^2=\|e_{x_0}-g_{\L_2,0}\|_+^2,
\end{equation}
for a unique $g_{\L_2,0}\in{\sf F}_{\text{loc},\L_2}$, where $\L_2:=R_2\cap\L$. Also, if we denote by
$\mcH_{2,+}$ the closure of $\cup_{\L\ssbs E}{\sf F}_{\text{loc},\L_2}$ w.r.t. the $\|\cdot\|_+$-norm, there is
a unique $g_{R_2,0}\in \mcH_{2,+}$ such that
\begin{equation}
\lim_{\L\ua E}g_{\L_2,0}=g_{R_2,0}\quad (\text{in }\mcH_+)
\end{equation}
and
\begin{equation}
\b=\|e_{x_0}-g_{R_2,0}\|_+^2.
\end{equation}
Similarly to \eqref{eq:orthogonality}, we have
\begin{equation}\label{eq:A^-1_operates}
A^{-1}(e_{x_0}-g_{R_2,0})=\b e_{x_0}+b_1\in \mcH_-,
\end{equation}
where $b_1$ is supported on $R_1$. Now we have on the one hand
\begin{eqnarray*}
_-\langle  e_{x_0}-f_{R_1,0},e_{x_0}-g_{R_2,0}\rangle_+&=&\lim_{\L\ua E} {}_-\langle  e_{x_0}-f_{\L_1,0},e_{x_0}-g_{\L_2,0}\rangle_+\\
&=&\lim_{\L\ua E} (e_{x_0}-f_{\L_1,0},e_{x_0}-g_{\L_2,0})_0\\
&=&\lim_{\L\ua E}1=1.
\end{eqnarray*}
On the other hand we have
\begin{eqnarray}\label{eq:VP_clue}
1&=&_-\langle  e_{x_0}-f_{R_1,0},e_{x_0}-g_{R_2,0}\rangle_+\nonumber\\
&=&{}_+\langle  A(e_{x_0}-f_{R_1,0}),A^{-1}(e_{x_0}-g_{R_2,0})\rangle_-\nonumber\\
&=&{}_+\langle  \a e_{x_0}+a_2,\b e_{x_0}+b_1\rangle_-\nonumber\\
&=&\a\b+{}_+\langle  a_2,b_1\rangle_-.
\end{eqnarray}
The proof is completed if we could show that ${}_+\langle
a_2,b_1\rangle_-=0$. Notice that $a_2$ is supported on $R_2$ and
$b_1$ on $R_1$, and $R_1\cap R_2=\es$. Thus it seems that
${}_+\langle a_2,b_1\rangle_-=0$, but we need to confirm it.

{\it Step 2: The case when $A$ is strictly positive}. Suppose that
there exist $c_1,\,c_2>0$ such that $c_1 I\le A\le c_2 I$. In this
case $A$ has a bounded inverse $A^{-1}$ in $\mcH_0=l^2(E)$. The RK
$B(x,y)$ for $\mcH_-$ (see \eqref{eq:kernel_B_recognization}) is
given by
\begin{equation}\label{eq:kernel_B_bounded_case}
B(x,y)={}_+\langle e_x,A^{-1}e_y\rangle_-=(e_x,A^{-1}e_y)_0,\quad
x,y\in E.
\end{equation}
Moreover, as for the elements, the inclusions in
\eqref{eq:inclusions} now become the equalities and the dual
pairings in \eqref{eq:dual_pair} and \eqref{eq:conjugate_form} are
just the inner product in the center space $\mcH_0$:
\begin{equation}
{}_-\langle  f,g\rangle_+=(f,g)_0=\ol{(g,f)_0}=\ol{{}_+\langle
g,f\rangle_-},\quad f,g\in \mcH_0=\mcH_-=\mcH_+,
\end{equation}
the equalities $\mcH_0=\mcH_-=\mcH_+$ meaning that all the spaces
have the same elements. We will, however, keep the pairing
notations ${}_-\langle  \cdot,\cdot\rangle_+$ and ${}_+\langle
\cdot,\cdot\rangle_-$ for a convenience. Now let us come back to
the equation \eqref{eq:VP_clue}. The dual pairing is just an inner
product in $\mcH_0$ and the vector $a_2$ vanishes on $R_1$ and
$b_1$ lives only on $R_1$. We therefore have
\begin{equation}\label{eq:exclusion_principle}
{}_+\langle  a_2,b_1\rangle_-=(a_2,b_1)_0=0.
\end{equation}
From \eqref{eq:VP_clue} and \eqref{eq:exclusion_principle} we have
$\a\b=1$.

We now extend the formula in Proposition
\ref{prop:results_of_local_systems}(b) to the infinite system.
That is, we will show that if $A$ is strictly positive, then
\begin{equation}\label{eq:alpha_realized}
\a=A(x_0,x_0)-A(x_0,R_1)A(R_1,R_1)^{-1}A(R_1,x_0).
\end{equation}
Notice that the function $A(\cdot,x_0)$ is an element of the space $\mcH_+$ and for each $\D\sbs E$ the function
$\D\ni y\mapsto A(y,x_0)$, which we denote by $A(\D,x_0)$, is the restriction of $A(\cdot, x_0)$ to the set
$\D$. Following Theorem \ref{thm:Aronszajn}, we denote this space by $\mcH_{\D;A}$ equipped with the norm
$\|\cdot\|_{\D;A}$. Since $A_\D(x,y)$, $x,y\in \D$, is the RK for $\mcH_{\D;A}$, it is obvious that for each
$g\in \mcH_{\D;A}$,
\begin{equation}
\|g\|_{\D;A}^2=(g,A(\D,\D)^{-1}g)_0,
\end{equation}
where $(\cdot,\cdot)_0$ is the usual inner product in $l^2(\D)$.
Thus \eqref{eq:alpha_realized} is equivalent to saying that
\begin{equation}\label{eq:alpha_realized2}
\a=A(x_0,x_0)-\|A(R_1,x_0)\|_{R_1;A}^2.
\end{equation}
Now by Proposition \ref{prop:results_of_local_systems}(b) we see
that for each $\L\ssbs E$, putting $\L_1=R_1\cap\L$,
\begin{eqnarray}\label{eq:alpha_local2}
\a_\L&=&A(x_0,x_0)-A(x_0,\L_1)A(\L_1,\L_1)^{-1}A(\L_1,x_0)\nonumber\\
&=&A(x_0,x_0)-\|A(\L_1,x_0)\|_{\L_1;A}^2.
\end{eqnarray}
On the other hand, as $\L$ increases, we have by Theorem \ref{thm:AronszajnI},
\begin{equation}\label{eq:convergence}
\lim_{\L\ua E}\|A(\L_1,x_0)\|_{\L_1;A}^2=\|A(R_1,x_0)\|_{R_1;A}^2.
\end{equation}
From \eqref{eq:alpha_and_beta} and
\eqref{eq:alpha_realized2}-\eqref{eq:convergence} we have shown
\eqref{eq:alpha_realized}.

{\it Step 3: The case when one of $R_1$ and $R_2$ is finite.} In this case either $a_2$ in \eqref{eq:A_operates}
or $b_1$ in \eqref{eq:A^-1_operates} is finitely supported. Moreover, since they have disjoint supports, by
\eqref{eq:dual_pair_realization} we have
\begin{equation}
{}_+\langle  a_2,b_1\rangle_-=\sum_{x\in E}\ol{a_2(x)}b_1(x)=0.
\end{equation}
This, together with \eqref{eq:VP_clue}, proves the theorem. This observation, as a matter of fact, gives us more
information. Notice that the number $\b$ in \eqref{eq:alpha_and_beta} is not altered even if we considered the
restriction of $B$ to the set $\wt{R_2}:=\{x_0\}\cup R_2$. Recall the notation $\|\cdot\|_{\wt{R_2};B}$ for the
norm in the RKHS $\mcH_{\wt{R_2};B}$ consisting of all the restrictions of vectors in $\mcH_-$ to the set
$\wt{R_2}$. $\mcH_{\wt{R_2};B}$ has its RK $B_{\wt{R_2}}(x,y)$, $x,y\in \wt{R_2}$, the restriction of $B$ onto
$\wt{R_2}$. We consider $\wt{R_2}$ being partitioned as $\wt{R_2}=\{e_{x_0}\}\cup \es\cup R_2$, and then apply
the result in this step to get
\begin{equation}\label{eq:beta_inverse}
\b^{-1}=\|e_{x_0}\|_{\wt{R_2};B}^2.
\end{equation}
In passing, we note that
$\|e_{x_0}\|_{\wt{R_2};B}^2=(e_{x_0},(B_{\wt{R_2}})^{-1}e_{x_0})_0$,
where $(B_{\wt{R_2}})^{-1}$ is the \lq\lq inverse\rq\rq of
$B_{\wt{R_2}}$ having the components
\begin{equation}
(B_{\wt{R_2}})^{-1}(x,y)=(e_x,e_y)_{\wt{R_2};B},\quad x,y\in
\wt{R_2}.
\end{equation}

For each $\e>0$, we introduce the strictly  positive and bounded operators $A(\e):=A+\e$ and $B(\e):=A(\e)^{-1}$
on $\mcH_0$. Let $\a(\e)$ and $\b(\e)$ be the numbers defined as in \eqref{eq:alpha_and_beta} by replacing the
operators $A$ and $B$ with $A(\e)$ and $B(\e)$, respectively. By the result in Step 2, we have
\begin{equation}\label{eq:VP_perturbation}
\a(\e)\b(\e)=1,\quad \e>0.
\end{equation}
On the other hand, by \eqref{eq:beta_inverse} we have
\begin{equation}
\b(\e)^{-1}=\|e_{x_0}\|_{\wt{R_2};B(\e)}^2:=(e_{x_0},(B(\e)_{\wt{R_2}})^{-1}e_{x_0})_0.
\end{equation}
In Lemma \ref{lem:limit_of_perturbation}, we have shown that
\begin{equation}\label{eq:limit_norm}
\lim_{\e\to
0}\|e_{x_0}\|_{\wt{R_2};B(\e)}^2=\|e_{x_0}\|_{\wt{R_2};B}^2,
\end{equation}
that is
\begin{equation}\label{eq:beta_inverse_limit}
\lim_{\e\to 0}\b(\e)^{-1}=\b^{-1}.
\end{equation}
It is easy to check that
\begin{equation}\label{eq:alpha_limit}
\lim_{\e\to 0}\a(\e)=\a.
\end{equation}
We thus get by \eqref{eq:VP_perturbation},
\eqref{eq:beta_inverse_limit}-\eqref{eq:alpha_limit}, $\a\b=1$.
The proof is completed. \EndProof

\section{Proofs of Theorem \ref{thm:global_PI} and Theorem \ref{thm:Gibbs_measure}}\label{sec:proof of global_PI}
The proof of Theorem \ref{thm:global_PI} will follow from the variational principle of Theorem
\ref{thm:variational_principle} and the projection-inversion inequalities, which we now introduce. For a matrix
$A$ on $E$, we denote by $A_\L$ for the submatrix, or projection of $A$ on the set $\L\sbs E$.
\begin{lem}\label{lem:monotonicity}
Let $A(x,y)$ and $B(x,y)$, $x,y\in E$, be the RK's respectively for $\mcH_+$ and $\mcH_-$ in Section 2. Then,
for any finite subsets $\L\sbs\D\ssbs E$, the  following inequalities hold:
\begin{enumerate}
\item[(a)] $(A_\L)^{-1}\le ((A_\D)^{-1})_\L\le B_\L$; \item[(b)] $(B_\L)^{-1}\le ((B_\D)^{-1})_\L\le A_\L$.
\end{enumerate}
\end{lem}
For a proof we need the projection-inversion lemma (see \cite[p 18]{OP}, \cite[Corollary 5.3]{ST2}, and
\cite[Lemma A.5]{GY}):
\begin{lem}\label{lem:projecton-inversion}
Let $T$ be any bounded positive definite operator with bounded inverse $T^{-1}$. Then for any projection $P$,
\begin{equation}\label{eq:projection-inversion}
P(PTP)^{-1}P\le PT^{-1}P.
\end{equation}
\end{lem}
\medskip
\Proof[ of Lemma \ref{lem:monotonicity}] The first inequalities in (a) and (b) follow from Lemma
\ref{lem:projecton-inversion}. In order to prove the second inequalities it is enough to show $(A_\L)^{-1}\le
B_\L$ and $(B_\L)^{-1}\le A_\L$, because $(B_\D)_\L=B_\L$ and $(A_\D)_\L=A_\L$ for $\L\sbs\D\ssbs E$. Moreover,
since the matrices are positive definite, either one of the inequalities $(A_\L)^{-1}\le B_\L$ or
$(B_\L)^{-1}\le A_\L$ implies the other. So, it is enough to prove $(B_\L)^{-1}\le A_\L$. Let $\|\cdot\|_{\L;B}$
be the norm on the space $\mcH_{\L;B}$ of all restrictions of functions of $\mcH_-$ to the subset $\L$ given in
Theorem \ref{thm:Aronszajn}. Since the function $B_\L(x,y)$, $x,y\in \L$, is the corresponding RK for
$\mcH_{\L;B}$, it is obvious that
\begin{equation}\label{eq:local_B_norm}
\|f_\L\|_{\L;B}^2=(f_\L,(B_\L)^{-1}f_\L)_0.
\end{equation}
On the other hand, since $\|f_\L\|_{\L;B}$ is the smallest number for all the values $\|g\|_-$ such that
$g_\L=f_\L$, we have the inequality
\begin{equation}\label{eq:bound_of_local_B_norm}
\|f_\L\|_{\L;B}^2\le \|f_\L\|_-^2=(f_\L,A_\L f_\L)_0.
\end{equation}
Combining \eqref{eq:local_B_norm} and \eqref{eq:bound_of_local_B_norm} we get the inequality $(B_\L)^{-1}\le
A_\L$, and the proof is completed. \EndProof
\begin{rem}\label{rem:extension}
Notice that $B$ is formally the inverse of $A$ (see
\eqref{eq:kernel_B_recognization}), and that the operator $A$ may
have $0$ in its spectrum (it should then be a continuous
spectrum). In that case, the operator $B$, considered on the space
$\mcH_0$, is an unbounded operator. Therefore, Lemma
\ref{lem:monotonicity} extends Lemma
\ref{lem:projecton-inversion}.
\end{rem}

Next we discuss the order relations between the restriction operators and the interaction operators giving the
local probability densities of DPP's. For each finite set $\L\sbs E$, we let, as before,
\begin{equation}
A_\L:=P_\L AP_\L\text{ and }A_{[\L]}:=K_\L(I-K_\L)^{-1},
\end{equation}
where $P_\L$ is the projection on $\mcH_0=l^2(E)$ onto $l^2(\L)$ and $K:=A(I+A)^{-1}$. We let
\begin{equation}
B_{[\L]}:=(A_{[\L]})^{-1}
\end{equation}
and recall that $B$ is the inverse of $A$.
\begin{lem}\label{lem:bounded_by_restriction}
For any finite set $\L\sbs E$,
\begin{equation}\label{eq:restriction_ordering}
A_{[\L]}\le A_\L \quad\text{and}\quad B_{[\L]}\le B_\L.
\end{equation}
\end{lem}
\medskip
\Proof We first prove the inequality $A_{[\L]}\le A_\L$. By Lemma \ref{lem:projecton-inversion},
\begin{eqnarray}\label{eq:A_Lambda_bound}
A_{[\L]}&=&-I_\L+((I-K)_\L)^{-1}\nonumber\\
&\le &-I_\L+P_\L(I-K)^{-1}P_\L\nonumber\\
&=&A_\L,
\end{eqnarray}
where $I_\L:=P_\L IP_\L$. The second inequality in
\eqref{eq:restriction_ordering} can be shown in two ways. We
introduce both of them. First, as before, we define $A(\e)=A+\e$
and $B(\e)=A(\e)^{-1}$ for $\e>0$. By the same way used in
\eqref{eq:A_Lambda_bound} we can show
\begin{equation}\label{eq:B_epsilon_bound}
B(\e)_{[\L]}:=(A(\e)_{[\L]})^{-1}\le B(\e)_\L,
\end{equation}
where $A(\e)_{[\L]}:=K(\e)_\L(I-K(\e)_\L)^{-1}$ with $K(\e):=A(\e)(I+A(\e))^{-1}$. Since $K(\e)\to K$ uniformly
as $\e\to 0$ we have
\begin{equation}
\lim_{\e\to 0}B(\e)_{[\L]}=(A_{[\L]})^{-1}=B_{[\L]}.
\end{equation}
On the other hand, by \eqref{eq:pointwise_convergence}
\begin{equation}\label{eq:restricted_B_epsilon_convergence}
B(\e)_\L\to B_\L \quad \text{uniformly as }\e\to 0.
\end{equation}
The inequality $B_{[\L]}\le B_\L$ follows from
\eqref{eq:B_epsilon_bound}-\eqref{eq:restricted_B_epsilon_convergence}.

The second way is to use Lemma \ref{lem:monotonicity}. $B_{[\L]}$ can be rewritten as
$B_{[\L]}=-I_\L+(K_\L)^{-1}$. Since $K=A(I+A)^{-1}$, $K$ satisfies the conditions in the hypothesis (H) of
Section 2. Applying Lemma \ref{lem:monotonicity}(a) for the pair of operators $K$ and $K^{-1}$, we have
\begin{eqnarray*}
B_{[\L]}&\le &-I_\L+P_\L K^{-1}P_\L\\
&=&P_\L(-I+K^{-1})P_\L\\
&=&P_\L A^{-1}P_\L=B_\L.
\end{eqnarray*}
The proof is completed. \EndProof \\We are now ready to prove Theorem \ref{thm:global_PI}.

\medskip
\Proof[ of Theorem \ref{thm:global_PI}] Let $x_0\in E$ and $\x\in \mcX$ be any configuration with $x_0\notin
\x$. For a convenience we define an auxiliary configuration $\ol \x\in \mcX$ as
\begin{equation}
\ol \x:=E\sm(\x\cup\{x_0\}).
\end{equation}
From the definition \eqref{eq:alpha_local_by_determinants} and Proposition
\ref{prop:results_of_local_systems}(b) we have the equality:
\begin{equation}
\a_{[\L]}=(A_{[\L]}(x_0\x_\L,x_0\x_\L)^{-1}(x_0,x_0))^{-1}.
\end{equation}
By the first inequality in Lemma \ref{lem:bounded_by_restriction}
and using Proposition \ref{prop:results_of_local_systems} once
more we have the bound
\begin{equation}\label{eq:bounded_by_alpha_local}
\a_{[\L]}\le (A_\L(x_0\x_\L,x_0\x_\L)^{-1}(x_0,x_0))^{-1}=\a_\L,
\end{equation}
where $\a_\L$ is defined in \eqref{eq:alpha_local} with $R_1:=\x$
(and $\L_1=\L\cap R_1=\L\cap \x\equiv \x_\L$). Now by Proposition
\ref{prop:results_of_local_systems}(b) and (c) we have
\begin{equation}
\b_{[\L]}:=(a_{[\L]})^{-1}=(B_{[\L]}(x_0\ol{\x}_\L,x_0\ol{\x}_\L)^{-1}(x_0,x_0))^{-1},
\end{equation}
where $B_{[\L]}=(A_{[\L]})^{-1}$. By the second inequality of
Lemma \ref{lem:bounded_by_restriction} we also have the bound
\begin{equation}\label{eq:bounded_by_beta_local}
\b_{[\L]}\le \b_\L,
\end{equation}
where, again, $\b_\L$ is defined in \eqref{eq:beta_local} with
$R_2:=\ol \x$. Now we take the limit of $\L$ increasing to the
whole space $E$. Since $\a_\L\to \a$ as $\L$ increases to $E$ we
have from \eqref{eq:bounded_by_alpha_local}
\begin{equation}\label{eq:limsup_bound}
\limsup_{\L\ua E}\a_{[\L]}\le \a.
\end{equation}
On the other hand, since $\b_\L\to \b$ as $\L\ua E$, we have also
from \eqref{eq:bounded_by_beta_local}
\begin{equation}\label{eq:liminf_bound}
\liminf_{\L\ua E}\a_{[\L]}=(\limsup_{\L\ua E}\b_{[\L]})^{-1}\ge
\b^{-1}=\a.
\end{equation}
The last equality comes from Theorem
\ref{thm:variational_principle}. From \eqref{eq:limsup_bound} and
\eqref{eq:liminf_bound} we get $\lim_{\L\ua E}\a_{[\L]}=\a$, which
was to be shown. \EndProof

Let us now turn to the proof of Theorem \ref{thm:Gibbs_measure}.
For the proof of Gibbsianness we will follow the method developed
in \cite{Y} for continuum models. We will first define a Gibbsian
specification \cite{G, Pr} by introducing an interaction. Then we
will prove that the DPP of our concern is admitted to the
specification. We refer also to \cite[Section 6]{ST2}. The proof
of uniqueness will be shown by following the method of \cite{ST2}.

Let $A$ be an operator that satisfies the conditions in the
hypothesis (H). For any finite configuration $\x\in \mcX$, we
define an interaction potential of the particles in $\x$ by
\cite{Y}
\begin{equation}\label{eq:potential_energy_again}
V(\x):=-\log\det A(\x,\x).
\end{equation}
Notice that $V(\x)>0$ for all finite configurations $\x\in \mcX$.
For any $\L_1,\L_2\ssbs E$ with $\L_1\cap \L_2=\emptyset$, and for
any configurations $\x_{\L_1}$ and $\x_{\L_2}$ on the sets $\L_1$
and $\L_2$, respectively, the mutual potential energy
$W(\x_{\L_1};\x_{\L_2})$ is defined to satisfy
\begin{equation}\label{eq:mutual_energy}
V(\x_{\L_1}\cup
\x_{\L_2})=V(\x_{\L_1})+V(\x_{\L_2})+W(\x_{\L_1};\x_{\L_2}).
\end{equation}
Now for each $\z_\L\in \mcX_\L$ and $\x\in \mcX$, we define the
energy of the particle configuration $\z_\L$ on $\L$ with boundary
condition $\x$ by
\begin{equation}\label{eq:potential_energy_with_b.c.}
H_\L(\z_\L;\x):=\lim_{\D\ua E}(V(\z_\L)+W(\z_\L;\x_{\D\sm\L})),
\end{equation}
whenever the limit exists. As a matter of fact, $H_\L(\z_\L;\x)$
is well-defined for all $\z_\L\in \mcX_\L$ and $\x\in \mcX$ as
shown in the following lemma:
\begin{lem}
Suppose that the operator $A$ satisfies the conditions in the
hypothesis (H). Then for any $\z_\L\in \mcX_\L$ and $\x\in \mcX$,
the value $H_\L(\z_\L;\x)$ in
\eqref{eq:potential_energy_with_b.c.} is well-defined as a finite
number.
\end{lem}
\Proof The proof is very similar to the one given for continuum
model in \cite[Lemma 3.2]{Y}. We define first for each bounded set
$\D\sps\L$
\begin{equation}
H_{\L;\D}(\z_\L;\x):=V(\z_\L)+W(\z_\L;\z_{\D\sm\L}).
\end{equation}
From the definitions
\eqref{eq:potential_energy_again}-\eqref{eq:mutual_energy} we get
\begin{equation}\label{eq:log_expression}
H_{\L;\D}(\z_\L;\x)=-\log\frac{\det
A(\z_\L\x_{\D\sm\L},\z_\L\x_{\D\sm\L})}{\det
A(\x_{\D\sm\L},\x_{\D\sm\L})}.
\end{equation}
Denoting $Q_\L$ for the projection on $l^2(\z_\L\x_{\L^c})$ onto
$l^2(\z_\L)$, $H_{\L;\D}(\z_\L;\x)$ can be rewritten as (cf.
Projection \ref{prop:results_of_local_systems})
\begin{equation}
H_{\L;\D}(\z_\L;\x)=-\log\det(Q_\L
A(\z_\L\x_{\D\sm\L},\z_\L\x_{\D\sm\L})^{-1}Q_\L )^{-1}.
\end{equation}
By using the projection-inversion lemma, Lemma
\ref{lem:projecton-inversion}, we see that $H_{\L;\D}(\z_\L;\x)$
decreases as $\D$ increases. Hence the limit
\begin{equation}\label{eq:local_energy_with_b.c.}
H_\L(\z_\L;\x)=\lim_{\D\ua E}H_{\L;\D}(\z_\L;\x)
\end{equation}
exists. We now show the finiteness of the limit value. Let
$\z_\L=\{x_1,\cdots,x_n\}$ be an enumeration of the sites in
$\z_\L$. Then we can rewrite the quantity inside the logarithm in
\eqref{eq:log_expression} as
\begin{eqnarray}\label{eq:determinant_decomposition}
&&\frac{\det A(\z_\L\x_{\D\sm\L},\z_\L\x_{\D\sm\L})}{\det
A(\x_{\D\sm\L},\x_{\D\sm\L})}\nonumber\\
&=&\frac{\det
A(x_1,\cdots,x_n\x_{\D\sm\L},x_1,\cdots,x_n\x_{\D\sm\L})}{\det
A(x_2,\cdots,x_n\x_{\D\sm\L},x_2,\cdots,x_n\x_{\D\sm\L})}
\cdots\frac{\det A(x_n\x_{\D\sm\L},x_n\x_{\D\sm\L})}{\det
A(\x_{\D\sm\L},\x_{\D\sm\L})}.
\end{eqnarray}
By Theorem \ref{thm:variational_principle}, each term in the
r.h.s. converges to a strictly positive number as $\D$ increases
to $E$. The proof is complete. \EndProof

The finiteness of the values $H_\L(\z_\L;\x)$ for all $\z_\L\in
\mcX_\L$ and $\x\in \mcX$ says that any configuration $\x\in \mcX$
is \lq\lq physically possible\rq\rq as noted in \cite[p 16]{Pr}.

Let us now define the Gibbsian specification. Define a partition
function on the set $\L$ with a boundary condition $\x\in \mcX$ as
\begin{equation}
Z_\L(\x):=\sum_{\z_\L\sbs \L}\exp[-H_\L(\z_\L;\x)].
\end{equation}
Then we define a probability distribution on the particle
configurations as
\begin{equation}\label{eq:specification_density}
\g_\L(\z_\L;\x):=\frac1{Z_\L(\x)}\exp[-H_\L(\z_\L;\x)].
\end{equation}

Let the set $\{0,1\}$ be equipped with a discrete topology and $\O:=\{0,1\}^E$ with a product topology. Let
$\mcF$ be the Borel $\s$-algebra on $\O$. For any subset $\D\sbs E$ we let $\mcF_\D$ be the $\s$-algebra on $\O$
such that the map $\x_{x}=1$ is measurable for all $x\in \D$. We notice that $\mcF_E=\mcF$. By the natural
mapping between $\O$ and $\mcX$, we define $\s$-algebras $\mcF_\D$, $\D\sbs E$, and $\mcF$ on $\mcX$. The
Gibbsian specification is defined as follows \cite{G,Pr}: for any measurable set $A\in \mcF$ and $\x\in \mcX$,
we define
\begin{equation}\label{eq:specification}
\g_\L(A|\x):=\sum_{\z_\L\sbs\L}\g_\L(\z_\L;\x)1_A(\z_\L\x_{\L^c}),
\end{equation}
where $1_A$ denotes the indicator function on the set $A$. It is not hard to check that the system
$(\g_\L)_{\L\ssbs E}$ defines a specification, i.e., it satisfies the following properties:
\begin{enumerate}
\item[(i)] $\g_\L(\cdot|\x)$ is a probability measure for each
$\x\in \mcX$; \item[(ii)] $\g_\L(A|\cdot)$ is
$\mcF_{\L^c}$-measurable for all $A\in \mcF$; \item[(iii)]
$\g_\L(A|\cdot)=1_A(\cdot)$ if $A\in \mcF_{\L^c}$; \item[(iv)]
$\g_\D\g_\L(A|\x):=\sum_{\z_\D\sbs \D}\g_\D(\z_\D|\x)\g_\L(A|\z_\D
\x_{\D^c})=\g_\D(A|\x)$ for all $\L\sbs\D\ssbs E$ and $\x\in
\mcX$.
\end{enumerate}
A probability measure $\m$ on $(\mcX,\mcF)$ is said to be {\it
admitted} to the specification $(\g_\L)_{\L\ssbs E}$, or a {\it
Gibbs measure}, if it satisfies the DLR-equations:
\begin{equation}
\m(A)=\int \g_\L(A|\x)d\m(\x),\quad\text{for any }A\in \mcF \text{
and }\L\ssbs E.
\end{equation}
The DLR condition says that for any $\L\ssbs E$ and $A\in \mcF$, the conditional expectation
$E^\m[1_A|\mcF_{\L^c}]$ has a version $\g_\L(A|\cdot)$:
\begin{equation}\label{eq:DLR_condition}
E^\m[1_A|\mcF_{\L^c}](\x)=\g_\L(A|\x),\quad \m\text{-a.a. }\x.
\end{equation}
From the equation \eqref{eq:determinant_decomposition} we easily
see that for any $\L\ssbs E$, $\z_\L\equiv\{x_1,\cdots,x_n\}\in
\mcX_\L$, and $\x\in \mcX$,
\begin{eqnarray*}
&&H_\L(\z_\L;\x)\\
&=&H_{\{x_1\}}(\{x_1\};\x)+H_{\{x_2\}}(\{x_2\};\{x_1\}\cup\x)
+H_{\{x_n\}}(\{x_n\};\{x_1,\cdots,x_{n-1}\}\cup\x).
\end{eqnarray*}
This says that all the values $H_\L(\z_\L;\x)$ are determined by
the values $H_{\{x\}}(\{x\};\x)$. Now then the DLR condition
\eqref{eq:DLR_condition} is equivalent to saying that (cf.
\cite{Sh} and \cite[Section 6]{ST2})
\begin{equation}\label{eq:Gibbs_property}
\frac{E^\m[\x_x=\{x\}|\mcF_{\{x\}^c}](\x)}{E^\m[\x_x=\es|\mcF_{\{x\}^c}](\x)}=\exp[-H_{\{x\}}(\{x\};\x)],
\quad \forall\,\, x\in E.
\end{equation}
\medskip
\Proof[ of Theorem \ref{thm:Gibbs_measure}] {\it Gibbsianness}. As noted above, it is enough to show the
relation \eqref{eq:Gibbs_property}. Let $x_0\in E$ be a fixed point and let $\x\in \mcX$. Then by
\eqref{eq:log_expression} and \eqref{eq:local_energy_with_b.c.},
\begin{equation}\label{eq:Boltzman_factor}
\exp[-H_{\{x_0\}}(\{x_0\};\x)]=\lim_{\D\ua E}\frac{\det A(x_0\x_{\D\sm \{x_0\}}, x_0\x_{\D\sm \{x_0\}})}{\det
A(\x_{\D\sm \{x_0\}}, \x_{\D\sm \{x_0\}})}.
\end{equation}
On the other hand, by
\eqref{eq:conditional_probability}-\eqref{eq:alpha_local_by_determinants}
\begin{eqnarray}\label{eq:PI}
\frac{E^\m[\x_{\{x_0\}}=\{x_0\}|\mcF_{\{x_0\}^c}](\x)}{E^\m[\x_{\{x_0\}}=\es|\mcF_{\{x_0\}^c}](\x)} &=&
\lim_{\D\ua E}\frac{E^\m[\x_{\{x_0\}}=\{x_0\}|\mcF_{\D\sm\{x_0\}}](\x_{\D\sm\{x_0\}})}
{E^\m[\x_{\{x_0\}}=\es|\mcF_{\D\sm\{x_0\}}](\x_{\D\sm\{x_0\}})}
\nonumber\\
&=&\lim_{\D\ua E}\frac{\det A_{[\D]}(x_0\x_{\D\sm \{x_0\}}, x_0\x_{\D\sm \{x_0\}})}{\det A_{[\D]}(\x_{\D\sm
\{x_0\}}, \x_{\D\sm \{x_0\}})}.
\end{eqnarray}
By Theorem \ref{thm:global_PI} the two limits in
\eqref{eq:Boltzman_factor} and \eqref{eq:PI} are the same and this
proves that the DPP $\m$ corresponding to the operator
$A(I+A)^{-1}$ is a Gibbs measure admitted to the specification
$(\g_\L)_{\L\ssbs E}$ in
\eqref{eq:specification_density}-\eqref{eq:specification}.

{\it Uniqueness}. Let us now address to the uniqueness problem of
the Gibbs measure. The arguments in the sequel parallel those in
\cite[Section 6]{ST2}. Suppose that $\n$ is a probability measure
admitted to the specification $(\g_\L)_{\L\ssbs E}$, i.e., $\n$
satisfies the condition \eqref{eq:DLR_condition}:
\begin{equation}\label{eq:nu_meets_DLR}
E^\n[1_A|\mcF_{\L^c}](\x)=\g_\L(A|\x),\quad \n\text{-a.a. }\x\in
\mcX \text{ for all }\L\ssbs E.
\end{equation}
Let $F:\mcX\to \mbR$ be a function of the form
\begin{equation}
F(\x)=1_{\{\x_{\L_0}=X\}},\quad \text{for some }\L_0\ssbs E \text{
and }X\sbs \L_0.
\end{equation}
We will show that for such functions $F$,
\begin{equation}\label{eq:nu=mu}
\n(F)=\m(F).
\end{equation}
Since those functions $F$ generate the $\s$-algebra $\mcF$, $\n$ then should be $\m$ and the uniqueness follows.

Let $\L\ssbs E$ be any set with $\L_0\sbs \L$. Then by \eqref{eq:nu_meets_DLR}
\begin{equation}\label{eq:CE_evaluated}
E^\n[F|\mcF_{\L^c}](\x)=\frac1{Z_\L(\x)}\sum_{Y\sbs\L\sm\L_0}\exp[-H_\L(X\cup
Y;\x)].
\end{equation}
Notice that the partition function $Z_\L(\x)$ can be rewritten as follows. Let $\Ph^{(\L;\x)}$ be a matrix of
size $|\L|$ whose components are given by
\begin{equation}
\Ph^{(\L;\x)}(x,y):=A(x,y)-(P_{\x_{\L^c}}A(\cdot,x),P_{\x_{\L^c}}A(\cdot,y)
)_{\x_{\L^c};A},
\end{equation}
where, as before, $A(\cdot,x)$ is a function on $E$:
$A(\cdot,x)(z)=A(z,x)$, $z\in E$, which belongs to $\mcH_+$, and
$P_{\x_{\L^c}}$ is the restriction operator restricting the
functions on $E$ to the set $\x_{\L^c}$, and
$(\cdot,\cdot)_{\x_{\L^c};A}$ is the inner product of the RKHS
$\mcH_{\x_{\L^c};A}$ with RK $A_{\x_{\L^c}}$, the restriction of
$A$ to the set $\x_{\L^c}$. By Theorem \ref{thm:Aronszajn}, the
matrix $ \Ph^{(\L;\x)}$ is well-defined. In an informal level, we
can write $\Ph^{(\L;\x)}(x,y)$ as
\begin{equation}
\Ph^{(\L;\x)}(x,y)=A(x,y)-A(x,\x_{\L^c})A(\x_{\L^c},\x_{\L^c})^{-1}A(\x_{\L^c},y).
\end{equation}
We refer to \cite[p 1559]{ST2} for the same matrix, where,
however, $A$ is strictly positive. For each finite $\D\sps \L$, we
let
\begin{eqnarray*}
\Ph^{(\L,\D;\x)}(x,y)&:=&A(x,y)-(P_{\x_{\D\sm\L}}A(\cdot,x),P_{\x_{\D\sm\L}}A(\cdot,y) )_{\x_{\D\sm\L};A},\\
&=&A(x,y)-A(x,\x_{\D\sm\L})A(\x_{\D\sm\L},\x_{\D\sm\L})^{-1}A(\x_{\D\sm\L},y),\quad x,y\in \L.
\end{eqnarray*}
By Theorem \ref{thm:AronszajnI},
\begin{equation}
\lim_{\D\ua E}\Ph^{(\L,\D;\x)}(x,y)=\Ph^{(\L;\x)}(x,y), \quad
x,y\in \L.
\end{equation}
Moreover, it is obvious that for any $X\sbs \L $
\begin{equation}
\exp[-H_{\L,\D}(X;\x)]=\det (\Ph^{(\L,\D;\x)}(X,X))
\end{equation}
and hence
\begin{equation}\label{eq:BF_by_determinant}
\exp[-H_{\L}(X;\x)]=\det (\Ph^{(\L;\x)}(X,X)).
\end{equation}
Therefore we get
\begin{equation}\label{eq:PF_reexpression}
Z_\L(\x)=\sum_{X\sbs \L}\exp[-H_\L(X;\x)]=\sum_{X\sbs \L}
\det(\Ph^{(\L;\x)}(X,X))=\det(I+\Ph^{(\L;\x)}).
\end{equation}
By using the expression \eqref{eq:BF_by_determinant} we see that
\begin{eqnarray}\label{eq:boundary_sum}
\sum_{Y\sbs \L\sm\L_0}\exp[-H_\L(X\cup Y;\x)]&=&\sum_{Y\sbs
\L\sm\L_0} \det(\Ph^{(\L;\x)}_{X\cup Y})\nonumber\\
&=&\det(P_{\L\sm\L_0}+\Ph^{(\L;\x)}_{X\cup(\L\sm\L_0)}).
\end{eqnarray}
Here we have put $\Ph^{(\L;\x)}_{X\cup Y}\equiv
\Ph^{(\L;\x)}(X\cup Y,X\cup Y)$, etc. We insert
\eqref{eq:PF_reexpression}-\eqref{eq:boundary_sum} into the r.h.s.
of \eqref{eq:CE_evaluated} and after a short computation we obtain
the expression for $E^\n[F|\mcF_{\L^c}](\x)$ in
\eqref{eq:CE_evaluated} (see \cite[eq. (6.47)]{ST2} for the
details):
\begin{eqnarray}\label{eq:CE_evaluation}
&&E^\n[F|\mcF_{\L^c}](\x)\nonumber\\
&=&\frac{\det(P_X[\Ph^{(\L;\x)}_{\L_0}-\Ph^{(\L;\x)}(\L_0,\L\sm\L_0)
(I+\Ph^{(\L;\x)}_{\L\sm\L_0})^{-1}\Ph^{(\L;\x)}(\L\sm\L_0,X)]P_X)}
{\det(I+\Ph^{(\L;\x)}_{\L_0}-\Ph^{(\L;\x)}(\L_0,\L\sm\L_0)
(I+\Ph^{(\L;\x)}_{\L\sm\L_0})^{-1}\Ph^{(\L;\x)}(\L\sm\L_0,\L_0))}.
\end{eqnarray}
We will show that
\begin{eqnarray}\label{eq:limit_of_denominator}
&&\lim_{\L\ua E}[
(I+\Ph^{(\L;\x)})_{\L_0}-\Ph^{(\L;\x)}(\L_0,\L\sm\L_0)
(I+\Ph^{(\L;\x)}_{\L\sm\L_0})^{-1}\Ph^{(\L;\x)}(\L\sm\L_0,\L_0)]\nonumber\\
&=& (I+A)_{\L_0}-A(\L_0,\L_0^c)
(I+A)(\L_0^c,\L_0^c)^{-1}A(\L_0^c,\L_0)\nonumber\\
&=&(P_{\L_0}(I+A)^{-1}P_{\L_0})^{-1}.
\end{eqnarray}
In fact, by using a similar computation as in Proposition
\ref{prop:results_of_local_systems}(b) we have for any $f_0\in
l^2(\L_0)$,
\begin{eqnarray}\label{eq:variational_computation}
&&(f_0,[(I+\Ph^{(\L;\x)})_{\L_0}-\Ph^{(\L;\x)}(\L_0,\L\sm\L_0)
(I+\Ph^{(\L;\x)}_{\L\sm\L_0})^{-1}\Ph^{(\L;\x)}(\L\sm\L_0,\L_0)]f_0)_{l^2(\L_0)}\nonumber\\
&=&\inf_{f\in
l^2(\L\sm\L_0)}(f_0-f,(I+\Ph^{(\L;\x)})(\L,\L)(f_0-f))_{l^2(\L)}\nonumber\\
&=& \inf_{f\in l^2(\L\sm\L_0)}(\|f_0-f\|_{l^2(\L)}^2+\inf_{g\in
l^2(\x_\L^c)}(f_0-f-g,A(f_0-f-g))_{l^2(E)})\nonumber\\
&=&\inf_{h\in l^2(\L_0^c)}(f_0-h,(P_\L+A)(f_0-h))_{l^2(E)}.
\end{eqnarray}
Since $P_\L\to I$ strongly as $\L\ua E$, it is obvious that the
last expression in \eqref{eq:variational_computation} converges as
$\L\ua E$ to
\begin{eqnarray}\label{eq:limit_of_v.c.}
&&\inf_{h\in
l^2(\L_0^c)}(f_0-h,(I+A)(f_0-h))_{l^2(E)}\nonumber\\
&=&(f_0,[(I+A)_{\L_0}-A(\L_0,\L_0^c)
(I+A)(\L_0^c,\L_0^c)^{-1}A(\L_0^c,\L_0)]f_0)_{l^2(\L_0)}.
\end{eqnarray}
Eqs. \eqref{eq:variational_computation}-\eqref{eq:limit_of_v.c.}
prove \eqref{eq:limit_of_denominator}. Recall the operator
$K=A(I+A)^{-1}$ which gives the DPP $\m$. We have
\begin{eqnarray}\label{eq:matrix_computation}
(I-K)_{\L_0}&=&P_{\L_0}(I+A)^{-1}P_{\L_0}\nonumber\\
&=&[(I+A)_{\L_0}-A(\L_0,\L_0^c)
(I+A)(\L_0^c,\L_0^c)^{-1}A(\L_0^c,\L_0)]^{-1}
\end{eqnarray}
and
\begin{eqnarray}\label{eq:PI_matrix}
A_{[\L_0]}=\frac{K_{\L_0}}{(I-K)_{\L_0}}&=&(I-K)_{\L_0}^{-1}\nonumber\\
&=&A(\L_0,\L_0)-A(\L_0,\L_0^c)
(I+A)(\L_0^c,\L_0^c)^{-1}A(\L_0^c,\L_0).
\end{eqnarray}
We thus get, by using \eqref{eq:CE_evaluation}-\eqref{eq:limit_of_denominator} and
\eqref{eq:matrix_computation}-\eqref{eq:PI_matrix},
\begin{eqnarray*}
\n(F)&=&\lim_{\L\ua E}E^\n[F|\mcF_{\L^c}](\x)\\
&=&\det (I-K_{\L_0})\det(P_XA_{[\L_0]}P_X)\\
&=&\det (P_XK_{\L_0}+P_{\L_0\sm X}(I-K_{\L_0}))\\
&=&\m(F).
\end{eqnarray*}
Now then $\n$ must be $\m$ and we have proven the uniqueness of
the Gibbs measure. \EndProof

\renewcommand{\thesection}{\Alph{section}}
\setcounter{section}{0}
\section{Appendix}
In this Appendix we discuss the hypothesis (H) in Section 2 by giving some examples. For simplicity we take
$E:=\mbZ$, the set of integers. We give three typical examples.

(i) {\it The case that $A$ is bounded and has a bounded inverse.} In this case all the norms $\|\cdot\|_-$,
$\|\cdot\|_0$, and $\|\cdot\|_+$ are equivalent and the spaces $\mcH_-$, $\mcH_0$, and $\mcH_+$ are the same as
sets. Obviously, $\mcH_-$ is functionally completed. The Gibbsianness of the DPP for the operator $A(I+A)^{-1}$
with $A$ being in this category has already been shown by Shirai and Takahashi \cite{ST2}.

(ii) {\it The case of diagonal matrices.} Suppose that $A$ is a diagonal matrix with diagonal elements $\a_x>0$
with $\a_x$ being bounded and decreasing to zero as $x\to \infty$. It is not hard to show that the hypothesis
(H) is satisfied for those operators $A$. In fact, $\mcH_-$ consists of those functions $f:E\to \mbC$ such that
$\sum_{x\in E}\a_x|f(x)|^2<\infty$. In other words, if $g=(g(x))_{x\in E}\in \mcH_0$ is any element of $\mcH_0$
then the vector $f\equiv (\a_x^{-1/2}g(x))_{x\in E}$ belongs to $\mcH_-$ and all the elements of $\mcH_-$ are of
this type.

(iii) {\it Perturbation of diagonal matrices.} Let $D$ be any
diagonal matrix of the type in the case (ii) above. Let
$A:=C^*DC$, where $C$ is a matrix such that $C$ and its inverse
$C^{-1}$ have off-diagonal elements that decrease sufficiently
fast as the distance from the diagonal become far. To say more
concretely, let $C(x,y)$ and $C^{-1}(x,y)$ be the matrix
components of $C$ and $C^{-1}$, respectively. We assume that there
exist positive numbers $m>0$ and $M>0$ such that
\begin{equation}
m\le C(x,x)\le M\text{ and }m\le C^{-1}(x,x)\le M \text{ for all
}x\in E,
\end{equation}
and $C(x,y)$ and $C^{-1}(x,y)$ converge to zero sufficiently fast
as $|x-y|\to \infty$. Then $A$ satisfies the conditions in (H).
Here we give an example. Let $D$ be a diagonal matrix with
diagonal elements $\a_x$, $x\in E$. We assume that there is $k\in
\mbN$ such that
\begin{equation}
\a_x^{-1}\le (1+|x|)^k,\quad x\in E.
\end{equation}
Let $C$ be a bounded operator with bounded inverse $C^{-1}$ such
that there is $m\ge 2(k+1)$ and
\begin{equation}\label{eq:off-diagonal_bound}
|C(x,y)|\le \frac1{1+|x-y|^m}\text{ and }|C^{-1}(x,y)|\le
\frac1{1+|x-y|^m}.
\end{equation}
Such an operator $C$ can, for example, be obtained by taking its convolution kernel function as the Fourier
series of strictly positive and sufficiently smooth function on the circle. We prove that $\mcH_-$ is
functionally completed. It is enough to show that the pre-Hilbert space $(\mcH_0,\|\cdot\|_-)$ satisfies the
conditions (i) and (ii) of Theorem \ref{thm:Aronszajn}. First we show that for any $y\in E$, $f(y)$ is
continuous in $(\mcH_0,\|\cdot\|)_-$. As noted in the Subsection 2.1, it is equivalent to show that $e_y\in
\mcH_+$. But, we have
\begin{eqnarray}\label{eq:continuity_of_the_e.m.}
\|e_y\|_+^2=(e_y,A^{-1}e_y)_0&=&((C^{-1})^*e_y,D^{-1}(C^{-1})^*e_y)_0\nonumber\\
&=&\sum_{x\in E}\a_x^{-1}|(C^{-1})^*e_y(x)|^2\nonumber\\
&\le &\sum_{x\in E}(1+|x|)^k\frac1{(1+|x-y|^m)^2}<\infty.
\end{eqnarray}
Next, notice that for any $f\in \mcH_0$,
\begin{equation}
\|f\|_-^2=(f,C^*DCf)_0=(Cf,DCf)_0=(\|Cf\|_-^{(D)})^2,
\end{equation}
where $\|\cdot\|_-^{(D)}$ is the \lq\lq $-$\rq\rq-norm for
$A\equiv D$. Since we have observed in case (ii) that the space
$\mcH_-^{(D)}$, completion of $\mcH_0$ w.r.t.
$\|\cdot\|_-^{(D)}$-norm, is functionally completed, it is enough
to show that given any sequence $\{f_n\}\sbs \mcH_0$ which is
$\|\cdot\|_-$-Cauchy and such that $f_n(y)\to 0$ as $n\to \infty$
for all $y\in E$, $Cf_n(y)\to 0$ as $n\to \infty$ for all $y\in
E$, because $\{Cf_n\}$ is $\|\cdot\|_-^{(D)}$-Cauchy. We observe
that there is a constant $b>0$ such that
\begin{equation}\label{eq:uniformly_bounded}
|f_n(y)|\le b(1+|y|)^k,\quad y\in E.
\end{equation}
In fact,
\begin{equation}
|f_n(y)|=|(e_y,f_n)_0|\le \|e_y\|_+\|f_n\|_-.
\end{equation}
Since $\{f_n\}$ is $\|\cdot\|_-$-Cauchy, $\|f_n\|_-$ is bounded
uniformly for $n\in \mbN$. On the other hand from
\eqref{eq:continuity_of_the_e.m.}, it is not hard to see that
there exists $b_1>0$ such that
\begin{equation}
\|e_y\|_+\le b_1(1+|y|)^k,\quad y\in E.
\end{equation}
This proves \eqref{eq:uniformly_bounded}. Now we have
\begin{eqnarray*}
Cf_n(y)&=&(e_y,Cf_n)_0\\
&=&(C^*e_y,f_n)_0\\
&=&\sum_{x\in E}\ol{C^*e_y(x)}f_n(x)\\
&=&\sum_{x\in \ptl(y)}\ol{C^*e_y(x)}f_n(x)+\sum_{x\in
\ptl(y)^c}\ol{C^*e_y(x)}f_n(x),
\end{eqnarray*}
where $\ptl(y)$ is any sufficiently large but finite set
containing $y$. Since $f_n(x)\to 0$ for all $x\in E$, the first
term in the last expression converges to $0$ as $n\to \infty$. By
using \eqref{eq:off-diagonal_bound} and
\eqref{eq:uniformly_bounded} we have
\begin{equation}\label{eq:small_limit}
|\sum_{x\in \ptl(y)^c}\ol{C^*e_y(x)}f_n(x)|\le b\sum_{x\in \ptl(y)^c}\frac1{1+|x-y|^m}(1+|x|)^k.
\end{equation}
Once $\ptl(y)$ has been taken sufficiently large, the quantity in
the r.h.s. of \eqref{eq:small_limit} becomes as much small as we
wish. This completes the proof.

\vskip 0.5 true cm
 \noindent{\it Acknowledgments.} The author thanks Prof. Y. M. Park and Dr. C. Bahn for fruitful discussions.

\end{document}